\documentclass[11pt]{article}

\usepackage{graphicx}
\usepackage{latexsym,amsmath,amsfonts,amscd, amsthm, dsfont}

\usepackage{bm,color}
\usepackage{epsfig,verbatim,epstopdf,graphics}
\usepackage{subfigure}
\usepackage{changebar}
\usepackage{multirow}

\usepackage{yhmath}%加弧线头标的字体包  \wideparen
 \usepackage{booktabs} %\cmidrule(lr){2-5} \cmidrule(lr){6-9}
 \usepackage{tikz}
\usepackage{verbatim}
\usetikzlibrary{arrows,backgrounds,snakes,shapes}
 \numberwithin{equation}{section}

\graphicspath{{./}{./figure/}}
\allowdisplaybreaks

\newtheoremstyle{plainNoItalics}{}{}{\normalfont}{}{\bfseries}{.}{ }{}

\theoremstyle{plain}
\newtheorem{thm}{Theorem}[section]

\theoremstyle{plainNoItalics}

\newtheorem{rem}[thm]{Remark}

\newtheorem{exa}[thm]{Example}

\newcommand{\ba}{{\bf a}}

\newcommand{\be}{\begin{eqnarray}}
\newcommand{\ee}{\end{eqnarray}}
\newcommand{\beno}{\begin{eqnarray*}}
\newcommand{\eeno}{\end{eqnarray*}}

\makeatletter

\newcommand{\Rmnum}[1]{\expandafter\@slowromancap\romannumeral #1@}
\makeatother

\begin{document}

\baselineskip=1.8pc

%\vspace*{.10in}

%=============  title  =========================
\begin{center}
{\bf
A high order conservative semi-Lagrangian discontinuous Galerkin method for two-dimensional transport simulations
}
\end{center}

\vspace{.2in}
\centerline{
Xiaofeng Cai\footnote{
 Department of Mathematics, University of Houston, Houston, TX, 77004. E-mail: xfcai@math.uh.edu.
},
 Wei Guo \footnote{
Department of Mathematics, Michigan State University, East Lansing, MI, 48824. E-mail:
wguo@math.msu.edu. Research
is supported by NSF grant NSF-DMS-1620047.
},
Jing-Mei Qiu\footnote{Department of Mathematics, University of Houston, Houston, TX, 77004. E-mail: jingqiu@math.uh.edu. Research supported by NSF grant NSF-DMS-1522777 and Air Force Office of Scientific Computing FA9550-12-0318.}
}

\bigskip
\centerline{
{\em Dedicated to Prof. Chi-Wang Shu on his 60th birthday.}
}
\bigskip
\noindent
{\bf Abstract.}
%\bigskip
In this paper, we develop a class of  high order conservative semi-Lagrangian (SL) discontinuous Galerkin (DG) methods for solving multi-dimensional linear transport equations. The  methods rely on a characteristic Galerkin weak formulation, leading to $L^2$ stable discretizations for linear problems. Unlike many existing SL methods,  the high order accuracy and mass conservation of the proposed methods are realized in a non-splitting manner. Thus, the detrimental splitting error, which is known to significantly contaminate long term transport simulations, will be not incurred. One key ingredient in the scheme formulation, borrowed from CSLAM {\em [Lauritzen, Nair \& Ullrich, 2010]}, is the use of Green's theorem which allows us to convert volume integrals into a set of line integrals. The resulting line integrals are much easier to approximate with high order accuracy, hence facilitating the implementation. Another novel ingredient is the construction of quadratic curves in approximating sides of upstream cell, leading to quadratic-curved quadrilateral upstream cells. Formal third order accuracy is obtained by such a construction. The desired positivity-preserving property is further attained by incorporating a high order bound-preserving filter.  To assess the performance of the proposed methods, we test and compare the numerical schemes with a variety of configurations for solving several benchmark transport problems with both smooth and nonsmooth solutions. The efficiency and efficacy are numerically verified.

\vfill

\newpage

\section{Introduction}

Transport phenomena are ubiquitous in nature, which may be described by a set of transport equations. We are concerned with the following first order transport equation
\begin{equation}
\label{eq:trans}
\frac{\partial u}{\partial t} + \nabla\cdot(\mathbf{a}u)  = 0,
\end{equation}
where $\mathbf{a}$ is the advection coefficient and could depend on space, time and even $u$ itself for a nonlinear problem. \eqref{eq:trans}  has a wide range of applications in science and engineering. For example, a real-world application of \eqref{eq:trans} is the multi-tracer transport process in the chemistry-climate model.  In such a model, the dynamical core (the fluid flow component) determines the wind field $\mathbf{a}$ that transports various physical and chemical substances in the atmosphere, ofter referred to as tracers, via a set of transport equations, see \cite{morgenstern2010review,lamarque2008simulated}. The tracers, on the other hand, provide feedback to the fluid flow through the parametrization process. 
%
 %for example, the Vlasov equation which describes fundamental dynamics of collisionless plasmas is in the form of \eqref{eq:trans}. 
 %The coefficient $\mathbf{a}$ representing the spatial advection and velocity acceleration of charged particles is  determined by the macroscopic quantities from solution $u$ in conjunction with field equations.
  As our initial effort to develop a class of genuinely high order and efficient transport schemes, we restrict our attention to the linear transport equations, meaning that $\mathbf{a}$ is independent of solution $u$. The extension to general transport equations including the Vlasov  equation in plasma physics
 and the multi-tracer transport model on the sphere will be addressed in our subsequent papers.

A transport scheme must feature several essential properties to qualify for practical applications. First, the transport equation  \eqref{eq:trans} may exhibit complex solution structures. For instance, 
%a solution of the Vlasov equation may develop thin filaments due to the phase mixing \cite{cheng2012study}.
in the multi-tracer transport process, the distributions of
tracers often develop rich structures in space such as clouds.
 Hence, it often requires that the scheme used is able to effectively resolve the fine-scale structures. Second, the equation \eqref{eq:trans} may conserve several physical quantities on the partial differential equations level, such as total mass, momentum, and energy. It is thus of paramount significance to conserve those physical invariants on the discrete level. Note that a violation of mass conservation will lead to large deviation and eventually crash long-term transport simulations \cite{Qiu_Shu2}. Lastly, a transport scheme should be efficient in terms of computational cost. The current generation global climate models include $\mathcal{O}(100)$ tracer species in order to adequately represent complex physical and chemical processes  \cite{lamarque2008simulated}. 
% Therefore, such a transport modeling requires a very careful balance between accuracy and computational complexity. In order to design effective schemes of practical value for the multi-tracer and other important transport models,
 To this end,
  we propose to develop a class of novel genuinely high order schemes, motivated by the following two facts:
(1) Many popular transport solvers are in the mesh-based Eulerian framework (see, e.g., \cite{harten1987uniformly,shu1988efficient,jiang1996efficient}), which are known to suffer from the inherent CFL time step restriction for stability. However, the time scales of the physics of interest may be well above the  time step restriction associated with the stability requirement \cite{white2011high}.
%  On the other hand, as a most remarkable property of the semi-Lagrangian (SL) approach, this type of scheme is able to attain desirable accuracy with extra large time step evolution, leading to great savings in computational cost.
A semi-Lagrangian (SL) approach avoids this issue while achieving desirable accuracy with time step restrictions only set by the physical quantities \cite{gucclu2014arbitrarily}, leading to great computational efficiency. For pioneering work on high order semi-Lagrangian schemes, we refer to several classical work on characteristic-Galerkin or Lagrangian-Galerkin methods\cite{morton1988stability, staniforth1991semi, giraldo1998lagrange}. 
(2) The widely recognized discontinuous Galerkin (DG) transport schemes \cite{reed1973triangularmesh,cockburn1989tvb,ayuso2009discontinuous,cockburn1989tvb2,cockburn1998runge,cockburn91s1,cockburn90s4,cockburn1998local,nair2005discontinuous,heath2012discontinuous} with excellent conservation properties  are very effective in resolving complex solution structures, thanks to the discontinuous nature of the approximation space. On the contrary, a continuous finite element method tends to introduce excessive numerical diffusion by the restrictive continuity requirement, resulting in smears of the solution or spurious oscillations \cite{heath2012discontinuous}.
%(3) SL approaches are multi-tracer efficient, since a substantial fraction of computational cost is independent of the number of tracers \cite{lamarque2008simulated}.
Therefore, the proposed work aims at incorporating  DG spatial discretization into the SL framework in a genuinely high order way in order to take advantage of both.  Many existing SLDG schemes are designed based on a dimensional splitting strategy due to its simplicity, which comes at the cost of a  splitting error, see \cite{James,qiu2011positivity,Guo2013discontinuous}. Such an error may be significant and hence greatly contaminating the solutions for long term transport simulations \cite{christlieb2014high}. On the other hand, the SLDG scheme proposed in \cite{restelli2006semi} is based on a flux form and free of splitting error, but subject to time step restriction, which degrades its computational efficiency to some extent.

 In our earlier work \cite{Guo2013discontinuous}, an SLDG weak formulation was formulated for one-dimensional (1D) transport equations based on a characteristic Galerkin weak formulation, see, e.g., \cite{morton1988stability,douglas1982numerical,celia1990eulerian,herrera1993eulerian}. Such a  method consists in breaking the upstream cell into several intersection subintervals and then evaluating the underlying integral on each subinterval via a high order quadrature rule. It is worth noting that
  the extension to  multi-dimensional cases without splitting is doable but very involved in implementation,
since the shape of an upstream cell may be irregular and an accurate numerical quadrature for arbitrary geometry in multi-dimensions is difficult to construct. Recently,  Lee \emph{et al.} proposed a non-splitting characteristic DG formulation for the two-dimensional (2D) transport equations in \cite{lee2016high}. Their method relies on approximating a upstream cell with a set of intersection polygons, followed by breaking each polygon into several triangles, and then generating quadrature points over each triangle. The scheme is second order accurate and the numerical results provided there are promising. However, as we mentioned, the extension to third or higher order accuracy is challenging. In particular, to attain the formal third order accuracy, quadratic curves are required to approximate boundaries of a upstream cell. Consequently, one needs to construct adequately accurate quadrature for triangles with a curved boundary, which is demanding in implementation.

In this paper, we seek to develop a non-splitting SLDG scheme that is unconditionally stable, genuinely high order accurate, mass conservative, and relatively easy to implement for multi-dimensional transport equations. We start with a reformulation of our previous 1D SLDG scheme \cite{Guo2013discontinuous}, and noting that the reformulated scheme does not rely on a numerical quadrature rule but the exact integration by means of fundamental theorem of calculus. The scheme has the potential to be generalized to multi-dimensional cases, since we are able to take advantage of the multi-dimensional generalization of the fundamental theorem of calculus, e.g., Green's theorem for the 2D case. By doing so, we indeed convert the area integrals into a set of line integrals which are much easier to evaluate. Such a technique has been used in a conservative semi-Lagrangian multi-tracer transport finite volume scheme (CSLAM) \cite{lauritzen2010conservative}. Unlike the formulation proposed in \cite{lee2016high}, the newly proposed scheme can be naturally extended to third order accuracy, by evaluating line integrals defined on quadratic curves as approximations to upstream cells boundaries. 

Note that  \eqref{eq:trans} preserves the maximum principle if $\ba$ is non-divergent. On the discrete level, it is highly desired that the numerical scheme used is able to preserve such a property, at least the positivity of the solution, since nonphysical negative values of tracers may trigger unrealistic processes in the climate modeling \cite{williamson2007evolution}. To address the issue, we employ a high order bound-preserving (BP) filter \cite{zhangshu2010} to ensure the desirable positivity preserving property, while retaining the original high order accuracy.

The rest of the paper is organized as follows. In Section 2,  we reformulate the 1D high order conservative SLDG scheme which can be extended to multi-dimensional cases. In Section 3, based on the reformulation, we develop a non-splitting conservative high order SLDG method for 2D transport equations. Some implementation details including the search algorithm are discussed.  We present several numerical results in Section 4 to benchmark the proposed scheme in terms of accuracy, efficiency, performance as well as the mass conservation property. We conclude the paper in Section 5.

\section{1D SLDG Scheme for Linear Transport}

We start with the following 1D linear transport equation
\begin{equation}
\frac{\partial u}{\partial t} + \frac{\partial }{\partial x} (a(x,t)u )=0,
\quad x\in[x_a,x_b],
\end{equation}
with a given initial condition and subject to proper boundary conditions.
%In this paper, we focus on the periodic boundary condition.
We assume $a(x,t)$ is continuous with respect to $x$ and $t$.
The formulation is designed in a very similar spirit to the scheme proposed in \cite{Guo2013discontinuous}, but with a slightly different implementation strategy for the ease of extension to 2D problems discussed in the next subsection.

In order to formulate the schemes, we start with a partition of computational domain $x_a=x_{\frac12}< x_{\frac32}<\cdots< x_{M+\frac12} =x_b$.
Let $I_j=[x_{j-\frac12}, x_{j+\frac12} ]$ denote an element of length $\Delta x=x_{j+\frac12}-x_{j-\frac12}$.% and its center $x_j$.
We define the finite dimensional approximation space, $V_h^k = \{ v_h:  v_h|_{I_j} \in P^k(I_j) \}$, where $P^k(I_j)$ denotes the set of polynomials of degree at most $k$  with $k\geq0$. Note that if $k=0$, the scheme formulated below reduces to a first order SL finite volume scheme.
%Then a local solution on $I_j$ can be expressed  as a polynomial of degree $k$,
%\begin{equation}
%x\in I_j:\ u_h^n(x,t) = \mathop{ \sum_{l=1}^{n_k} } \widehat{u}_j^{(l)}(t) \psi_l(x),
%%=\mathop{ \sum_{i=1}^{K} u_h^n(x_n^i,t) l_n^i(x) }.
%\end{equation}
%where $n_k= k +1$.
%We adopt the scaled Legendre polynomials, $\{\psi_l(x)\}$:
%\begin{equation}
%\psi_1(x) = 1 ,\  \psi_2(x) = \frac{x-x_j}{\Delta x},\ \psi_3(x) = \left( \frac{x-x_j}{\Delta x} \right)^2-\frac{1}{12},\ \psi_4(x) = \left( \frac{x-x_j}{\Delta x} \right)^3 - \frac{3}{20}\left(\frac{x-x_j}{\Delta x}\right),\cdots.
%\end{equation}
%

We update the solution at the time level $t^{n+1}$ over a cell $I_j$ from the solution at $t^n$. We employ the weak formulation of characteristic Galerkin method proposed in \cite{childs1990characteristic}. Specifically, we let the time-dependent test function $\psi(x,t)$ satisfy the adjoint problem with $\forall \Psi \in P^k(I_j)$,
\begin{equation}
\begin{cases}
\psi_t + a(x,t) \psi_x =0 ,\\
\psi(t=t^{n+1}) = \Psi(x), \ t\in[t^n,t^{n+1} ],
\end{cases}
\label{final-value}
\end{equation}
where $t^n$ denotes the $n$-th time level. Let $\Delta t =t^{n+1}-t^n$ denote the time step.
We  remark that for the above advective form of equation, the solution stays constant along a trajectory. It was shown in \cite{Guo2013discontinuous} that
\begin{equation}
\frac{d}{dt}\int_{I_j(t)} u(x,t)\psi(x,t) dx =0,
\label{property}
\end{equation}
where $I_j(t)$ is a dynamic interval bounded by characteristics emanating form cell boundaries of $I_j$ at $t=t^{n+1}$.
An SL time discretization of \eqref{property} leads to
\begin{equation}
\int_{I_j} u^{n+1} \Psi dx = \int_{I_{j}^\star } u(x,t^n) \psi(x,t^n)dx,
\label{integral1}
\end{equation}
where $I_{j}^\star  = [x_{j-\frac12}^\star , x_{j+\frac12}^\star]$ with $x_{j\pm\frac12}^\star$ being the foots of trajectory emanating from
$(x_{j\pm\frac12}^\star,t^{n+1})$  at the time $t^n$.
In particular, to update the numerical solution $u^{n+1}$, the following procedure is performed.
\begin{description}
  \item[(1)] We choose $k+1$ interpolation points $x_{j,q}$, such as the Gauss-Lobatto points (the midpoint for $k=0$) over the interval $I_j$ and locate the feet $x_{j,q}^\star$ (see Figure \ref{schematic_1d_1}) , by numerically solving  the following final-value problem (trajectory equation):
  \begin{equation}
  \frac{d x(t) }{dt}  = a( x(t),t ) \ \text{with} \ x(t^{n+1}) = x_{j,q}
  \end{equation}
   with a high order numerical integrator such as a fifth order Runge-Kutta method given in \cite{butcher2008numerical}.
   % In our implementation, we use  a fifth-order Runge-Kutta method  which is given in the Appendix \ref{RK5}.
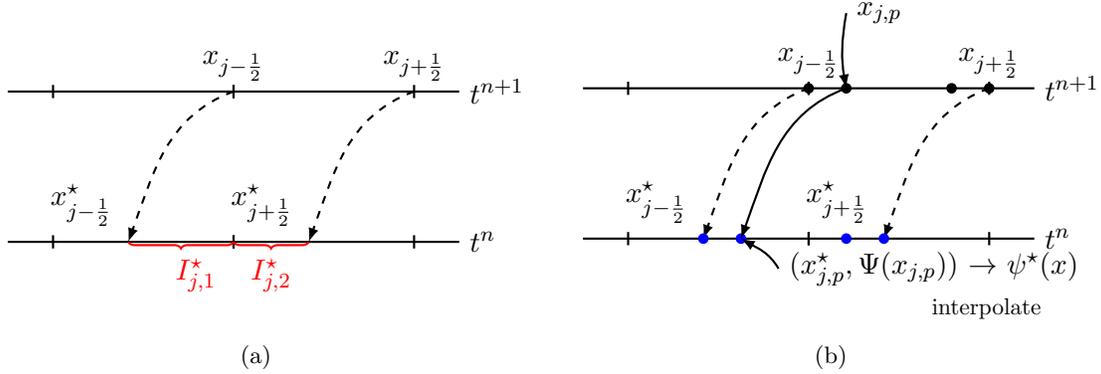
\begin{figure}[h!]
\centering
\subfigure[]{
\begin{tikzpicture}[x=1cm,y=1cm]%[scale=1.0]
  \begin{scope}[thick]
   \draw (-3,3) node[fill=white] {};
    \draw (-3,-1) node[fill=white] {};
    %\draw (60:-1cm) node[fill=white] {$E$} -- (60:3cm) node[fill=white] {$F$};
    \draw[black]                   (-3,0) node[left] {} -- (3,0)
                                        node[right]{$t^{n}$};
    \draw[black] (-3,2) node[left] {$$} -- (3,2)
                                        node[right]{$t^{n+1}$};

     %\draw[snake=ticks,segment length=2.2cm] (-2.4,0) -- (2.8,0);
     \draw[snake=ticks,segment length=2.4cm] (-2.4,2) -- (0,2) node[above] {$x_{j-\frac12}$};
     \draw[snake=ticks,segment length=2.4cm] (0,2) -- (2.4,2) node[above] {$x_{j+\frac12}$};

          \draw[snake=ticks,segment length=2.4cm] (-2.4,0) -- (0,0);
     \draw[snake=ticks,segment length=2.4cm] (0,0) -- (2.4,0);

     \draw[-latex,dashed](0,2)node[left,scale=1.3]{$$}
        to[out=200,in=70] (-1.4,0) node[above left=2pt] {$x_{j-\frac12}^\star$};  %%%thick 实线； dashed 虚线
             \draw[-latex,dashed](2.4,2)node[left,scale=1.3]{$$}
        to[out=200,in=70] ( 1.,0) node[above left=2pt] {$x_{j+\frac12}^\star$};

        \draw[snake=brace,mirror snake,red,thick] (-1.4,0) -- (0,0) node[below left=2.5pt] {$I_{j,1}^\star$};
        \draw[snake=brace,mirror snake,red,thick] (0,0) -- (1,0) node[below left=2.5pt] {$I_{j,2}^\star$};

  \end{scope}
\end{tikzpicture}
}
\subfigure[]{
\begin{tikzpicture}[x=1cm,y=1cm]%[scale=1.0]
  \begin{scope}[thick]
    %\draw (60:-1cm) node[fill=white] {$E$} -- (60:3cm) node[fill=white] {$F$};
    \draw[black]                   (0+1,0) node[left] {} -- (6+1,0)
                                        node[right]{$t^{n}$};
    \draw[black] (0+1,2) node[left] {$$} -- (6+1,2)
                                        node[right]{$t^{n+1}$};

     %\draw[snake=ticks,segment length=2.2cm] (-2.4,0) -- (2.8,0);
     \draw[snake=ticks,segment length=2.4cm] (0.6+1,2) -- (3+1,2) node[above] {$x_{j-\frac12}$};
     \draw[snake=ticks,segment length=2.4cm] (3+1,2) -- (5.4+1,2) node[above] {$x_{j+\frac12}$};

          \draw[snake=ticks,segment length=2.4cm] (0.6+1,0) -- (3+1,0);
     \draw[snake=ticks,segment length=2.4cm] (3+1,0) -- (5.4+1,0);

     \draw[-latex,dashed](3+1,2)node[left,scale=1.3]{$$}
        to[out=200,in=70] (1.6+1,0) node[above left=2pt] {$x_{j-\frac12}^\star$};  %%%thick 实线； dashed 虚线
             \draw[-latex,dashed](5.4+1,2)node[left,scale=1.3]{$$}
        to[out=200,in=70] ( 4+1,0) node[above left=2pt] {$x_{j+\frac12}^\star$};

%%%
\fill [black] (3+1,2) circle (2pt) node[] {};
\fill [black] (5.4+1,2) circle (2pt) node[] {};
\fill [black] (3+1+0.5,2) circle (2pt) node[] { };
\fill [black] (5.4+1-0.5,2) circle (2pt) node[] {};

\fill [blue] (3+1 -1.4,0) circle (2pt) node[] {};
\fill [blue] (5.4+1-1.4,0) circle (2pt) node[] {};
\fill [blue] (3+1+0.5-1.4,0) circle (2pt) node[] { };
\fill [blue] (5.4+1-0.5-1.4,0) circle (2pt) node[] {};

\draw[-latex,thick](3+1+0.5,3)node[right,scale=1.0]{$x_{j,p}$}
        to[out=260,in=100] (3+1+0.5,2) node[] {};

 \draw[-latex,thick](3+1+0.5,2)node[right,scale=1.0]{}
        to[out=200,in=70] (3+1+0.5-1.4,0) node[] {};

   \draw[-latex,thick](3+1+1-1.4,-0.4)node[right=0pt,scale=1.0]{$(x_{j,p}^\star,\Psi(x_{j,p}) )$ $\rightarrow$ $\psi^\star(x)$ }
   node[below right=8pt,scale=0.8]{$\ \ \ \ \ \ \ \ \ \ \ \ \ \ \ $ interpolate}
        to[out=120,in=330] (3+1+0.5-1.4,0) node[] {};

  \end{scope}
\end{tikzpicture}
}
\label{schematic_1d_1}
\caption{Schematic illustration for one-dimensional SLDG schemes.}
\end{figure}

\item[(2)] Recall that the test function $\psi$ solves the final-value problem \eqref{final-value} and hence stays constant along the characteristics, i.e., $\psi(x_{j,q}^\star ,t^n ) = \Psi(x_{j,q})$.
Now we are able to determine the unique polynomial $\psi^\star(x)$ of degree $k$ that interpolates $\psi(x,t^n)$, i.e. the test function at $t^n$, with the data points $(x_{j,q}^\star , \Psi(x_{j,q}))$,  $q=0,\cdots,k$, see Figure \ref{schematic_1d_1} (b).

  \item[(3)] Detect intervals/sub-intervals within $I_j^\star = \bigcup_l I_{j,l}^\star$, which are  the intersections between $I_j^\star$ and the grid elements ($l$ is the index for sub-interval). For instance, in Figure \ref{schematic_1d_1} (a), there are two sub-intervals: $I_{j,1}^\star = [x_{j-\frac12}^\star , x_{j-\frac12}] $
     and $I_{j,2}^\star = [x_{j-\frac12},x_{j+\frac12}^\star ]$.

\item[(4)] Lastly, the right hand side (RHS) of \eqref{integral1} is approximated by
\begin{equation}
 \int_{I_{j}^\star } u(x,t^n) \psi(x,t^n)dx \approx
\sum_l \int_{I_{j,l} } u^n(x) \psi^\star(x)dx.
\label{integral2}
\end{equation}
The summation is defined by incorporating all the sub-intervals of $I_j^\star$.
Note that the integrands in \eqref{integral2} are polynomials of degree $2k$; thus the integration can be evaluated exactly.

 \end{description}

The main difference between the above reinterpretation and the formulation proposed in \cite{Guo2013discontinuous} is that, rather than using the numerical quadrature to evaluate $ \int_{I_{j}^\star } u(x,t^n) \psi(x,t^n)dx $, we first
approximate $\psi(x,t^n)$ by an interpolating polynomial $\psi^\star(x)$, then evaluate the integrals exactly. Similarly, the newly proposed method is locally conservative (which can be easily verified by letting $\Psi(x)=1$). More importantly, based on the same idea, we can develop a non-splitting SLDG scheme for multi-dimensional transport equations by taking advantage of the multi-dimensional generalization of the fundamental theorem of calculus, e.g., Green's theorem in 2D cases.

%%%%%%%%%%%%%%%%%%%%%%%%%%%%%%%%%%%%%%%%%%%%%%%%%%%%%%%%%%%%%%%%%%%%%%%%%%%

\section{A Non-splitting SLDG Formulation for 2D Transport Problems  }
In this section, we develop a non-splitting SLDG scheme for 2D transport equations which is a natural generalization of the 1D reformulation we proposed in the previous section.
% We present SLDG  schemes with quadrilateral upstream cell and a third order SLDG scheme with quadratic-curved quadrilateral upstream cell for 2D problems. In the last, a bound-preserving filter is presented to enforce the positivity of the solutions.

%Next, we extend the 1D SLDG formulation to 2D problems.
Consider the following 2D transport problem
\begin{equation}
\frac{\partial u}{\partial t} + \frac{\partial}{\partial x} (a(x,y,t)u )
+ \frac{\partial }{\partial y} ( b(x,y,t)u ) = 0.
\label{2d_linear}
\end{equation}
%with a given initial condition.
Here $(a,b)$ is a prescribed   velocity field depending on time and space.
%For simplicity, we consider the periodic conditions.
We assume  a Cartesian partition of the computational domain $\Omega=\{ A_j \}_{j=1}^J$ (see Figure \ref{schematic_2d} (a)) for simplicity. Note that the procedure below also applies to unstructured meshes, but the implementation is more involved.
Similar to the 1D case, we define the DG approximation space as a finite dimensional vector space $V_h^k=\{ v_h:v_h|_{A_j} \in P^k(A_j) \}$, where $P^k(A_j)$ denotes the set of polynomials of degree at most $k$.
To update the numerical solution from time step $t^n$ to $t^{n+1}$ over cell $A_j$, we consider the following adjoint problem for the test function $\Psi\in P^k(A_j)$:
\begin{equation}
\psi_t + a(x,y,t) \psi_x + b(x,y,t) \psi_y =0,\
\text{subject to} \
\psi(t=t^{n+1}) = \Psi(x,y),\
t\in [t^{n},t^{n+1}].
\label{2d_adjoint}
\end{equation}

\begin{figure}[h]
\centering
\subfigure[]{
\begin{tikzpicture}
    %%%%%%%%%%%%%%%%%%%%%%%%%%%%%%%%%%%%%%%%%%%%%%%%%%%%%%%%%%mesh
    \draw[black,thin] (0,0.5) node[left] {} -- (5.5,0.5)
                                        node[right]{};
    \draw[black,thin] (0,2.) node[left] {$$} -- (5.5,2)
                                        node[right]{};
    \draw[black,thin] (0,3.5) node[left] {$$} -- (5.5,3.5)
                                        node[right]{};
    \draw[black,thin] (0,5 ) node[left] {$$} -- (5.5,5)
                                        node[right]{};
                                        %%%%%%%%%%%
    \draw[black,thin] (0.5,0) node[left] {} -- (0.5,5.5)
                                        node[right]{};
    \draw[black,thin] (2,0) node[left] {$$} -- (2,5.5)
                                        node[right]{};
    \draw[black,thin] (3.5,0) node[left] {$$} -- (3.5,5.5)
                                        node[right]{};
    \draw[black,thin] (5,0) node[left] {$$} -- (5,5.5)
                                        node[right]{};
    %%%%%%%%%%%%%%%%%%%%%%%%%%%%%%%%%%%%%%%%%%% background element
    \fill [blue] (3.5,3.5) circle (2pt) node[] {};
    \fill [blue] (5,3.5) circle (2pt) node[] {};
    \fill [blue] (3.5,5) circle (2pt) node[below right] {$A_j$} node[above left] {$v_4$};
    \fill [blue] (5,5) circle (2pt) node[] {};

     \draw[thick,blue] (3.5,3.5) node[left] {} -- (3.5,5)
                                        node[right]{};
      \draw[thick,blue] (3.5,3.5) node[left] {} -- (5,3.5)
                                        node[right]{};
       \draw[thick,blue] (3.5,5) node[left] {} -- (5,5)
                                        node[right]{};
        \draw[thick,blue] (5,3.5) node[left] {} -- (5,5)
                                        node[right]{};
    %%%%%%%%%%%%%%%%%%%%%%%%%%%%%%%%%%%%%%%%%%%%%%%%%%%
    %%%%%%%%%%%%%%%%%%%%%%%%%%%%%%%%%%%%%%%%%%%%%%%%%%% Lagrangian element
    \fill [red] (1.,1) circle (2pt) node[above right,black] {};
    \fill [red] (3,1) circle (2pt) node[] {};
    \fill [red] (1,2.5) circle (2pt) node[below right] {$A_j^\star$} node[above left] {$v_4^\star$};
    \fill [red] (2.5,2.5) circle (2pt) node[] {};

     \draw[-latex,dashed](3.5,5)node[right,scale=1.0]{}
        to[out=240,in=70] (1,2.50) node[] {};

     \draw (0.5+0.01,2-0.01) node[fill=white,below right] {$A_l$};

     \draw [red,thick] (1,1)node[right,scale=1.0]{}
        to[out=20,in=150] (2,0.7) node[] {};

        \draw [red,thick] (2,0.7)node[right,scale=1.0]{}
        to[out=330,in=240] (3,1) node[] {};
%********************************************************************
             \draw [red,thick] (1,2.5)node[right,scale=1.0]{}
        to[out=310,in=90] (1.1,2) node[] {};
        \draw [red,thick] (1.1,2)node[right,scale=1.0]{}
        to[out=270,in=80] (1,1) node[] {};

        \draw [red,thick] (1,2.5)node[right,scale=1.0]{}
        to[out=10,in=180] (2.5,2.5) node[] {};

        \draw [red,thick] (3,1)node[right,scale=1.0]{}
        to[out=80,in=280] (2.5,2.5) node[] {};
    %%%%%%%%%%%%%%%%%%%%%%%%%%%%%%%%%%%%%%%%%%%%%%%%%%%%%%%%%%%%%
\end{tikzpicture}
}
\subfigure[]{

\begin{tikzpicture}[scale = 1.3]
    %%%%%%%%%%%%%%%%%%%%%%%%%%%%%%%%%%%%%%%%%%%%%%%%%%%%%%%%%%mesh
    \draw[black,thin] (0,0.5) node[left] {} -- (4,0.5)
                                        node[right]{};
    \draw[black,thin] (0,2.) node[left] {$$} -- (4,2)
                                        node[right]{};
    \draw[black,thin] (0,3.5) node[left] {$$} -- (4,3.5)
                                        node[right]{};
                                        %%%%%%%%%%%
    \draw[black,thin] (0.5,0) node[left] {} -- (0.5,4)
                                        node[right]{};
    \draw[black,thin] (2,0) node[left] {$$} -- (2,4)
                                        node[right]{};
    \draw[black,thin] (3.5,0) node[left] {$$} -- (3.5,4)
                                        node[right]{};

    %%%%%%%%%%%%%%%%%%%%%%%%%%%%%%%%%%%%%%%%%%%%%%%%%%%
    %%%%%%%%%%%%%%%%%%%%%%%%%%%%%%%%%%%%%%%%%%%%%%%%%%% Lagrangian element
    \fill [red] (1.,1) circle (2pt) node[above right,black] {$A_{j,l}^\star$};
    \fill [red] (3,1) circle (2pt) node[] {};
    \fill [red] (1,2.5) circle (2pt) node[below right] {$A_j^\star$} node[above left] {};
    \fill [red] (2.5,2.5) circle (2pt) node[] {};

     \draw (0.5+0.01,2-0.01) node[fill=white,below right] {};

\draw [red,thick] (1,1)node[right,scale=1.0]{}
        to[out=20,in=150] (2,0.7) node[] {};

        \draw [red,thick] (2,0.7)node[right,scale=1.0]{}
        to[out=330,in=240] (3,1) node[] {};
%********************************************************************
             \draw [red,thick] (1,2.5)node[right,scale=1.0]{}
        to[out=310,in=90] (1.1,2) node[] {};
        \draw [red,thick] (1.1,2)node[right,scale=1.0]{}
        to[out=270,in=80] (1,1) node[] {};

        \draw [red,thick] (1,2.5)node[right,scale=1.0]{}
        to[out=10,in=180] (2.5,2.5) node[] {};

        \draw [red,thick] (3,1)node[right,scale=1.0]{}
        to[out=80,in=280] (2.5,2.5) node[] {};
        %%%%%%%%%%
           \draw (0.5+0.01,2-0.01) node[fill=white,below right] {$A_l$};
    %%%%%%%%%%%%%%%%%%%%%%%%%%%%%%%%%%%%%%%%%%%%%%%%%%%%%%%%%%%%%
         \draw[-latex,ultra thick] (1,1)node[right,scale=1.0]{}
        to  (2,1) node[] {};

             \draw[-latex,ultra thick]  (2,1)node[right,scale=1.0]{}
        to (2,2) node[] {};
             \draw[-latex,ultra thick]  (2,2)node[right,scale=1.0]{}
        to (1,2) node[] {};
             \draw[-latex,ultra thick]   (1,2)node[right,scale=1.0]{}
        to (1,1) node[] {};

            %%%%%%%%%%%%%%%%%%%%%%%%%%%%%%%%%%%%%%%%%%%%%%%%%%%%%%%%%%%%%
\draw [thick] (1,1)-- (3,1) node[] {};
 \draw [thick] (1,2.5) -- (1,1) node[] {};
        \draw [thick] (1,2.5)--(2.5,2.5) node[] {};
        \draw [thick] (3,1)--(2.5,2.5) node[] {};
\end{tikzpicture}

}
\caption{Schematic illustration of the SLDG formulation in two dimension. $P^1$ case.}
\label{schematic_2d}
\end{figure}

%
%We approximate $\Omega$ by $I\times J$ nonoverlapping elements, $(x,y)\in [x_i^l,x_i^r]\times[y_j^l,y_j^r]=A_{ij}$, $x_i = \frac{x_i^l + x_i^r}{2}$, $y_j = \frac{y_j^l + y_j^r }{2}$, $\Delta x = x_i^r - x_i^l$, $\Delta y = y_j^r - y_j^l$.
%
%
%We assume that the local solution in Eulerian cell $A_{ij}$ is expressed as a polynomial of degree $K$,
%\begin{equation}
%\ u_h^{ij} (x,y,t) = \sum_{k=1}^{n_K} \widehat{u}_{ij}^{(k)}(t) \psi_k (x,y,t),
%\end{equation}
%where $n_K = \frac{(K+1)(K+2)}{2}$.
%We adopt the scaled Legendre polynomials: $\psi_1(x,y) =1$, $\psi_2(x,y)=\frac{x-x_i}{\Delta x}$, $\psi_3(x,y) = \frac{y-y_j}{\Delta y}$, $\psi_4(x,y) = \left( \frac{x-x_i}{\Delta x} \right)^2 -\frac{1}{12}$, $\psi_5(x,y)=\frac{(x-x_i)(y-y_j) }{\Delta x\Delta y}$, $\psi_6(x,y) = \left( \frac{y-y_j}{\Delta y} \right)^2 -\frac{1}{12}$, $\cdots$.

The scheme formulation takes advantage of the identity
\begin{equation}
\frac{d}{dt} \int_{A_{j}(t) } u(x,y,t) \psi(x,y,t) dxdy =0,
\label{2d_dt}
\end{equation}
where $A_{j}(t)$ is a dynamic moving cell, emanating from the Eulerian cell $A_{j}$ at $t^{n+1}$ backward in time by following characteristics trajectories.
The non-splitting SLDG scheme is formulated as follows:
Given the approximate solution $u^n\in V_h^k$ at time $t^n$, find $u^{n+1}\in V_h^k$ such that $\forall \Psi \in V_h^k$, we have
\begin{equation}
\int_{A_j}u^{n+1} \Psi(x,y) dxdy =
\int_{A_j^\star} u^n \psi(x,y,t^n) dxdy,
\label{2d_temporal}
\end{equation}
for $j= 1,\cdots,J$ and $n=0,1,\cdots.$. $A_j^\star$ denotes the upstream cell of grid cell $A_j$ following the characteristics backward to $t^n$, see the deformed cell bounded by red curves in Figure \ref{schematic_2d} (a), and the test function $\psi(x, y, t)$ satisfies the adjoint problem \eqref{2d_adjoint}.

A key step of the proposed methodology is the evaluation of the RHS of \eqref{2d_temporal}.
 In general, the upstream cell of the RHS of \eqref{2d_temporal} is no longer a rectangle or a quadrilateral. In the following subsections, we will discuss several techniques to evaluate the volume integral in \eqref{2d_temporal}, where the shape of the upstream cell is approximated by either a quadrilateral or a quadratic-curved quadrilateral with the goal to achieve second or third order spatial accuracy.
%The procedure to evaluate the volume integral on the RHS of \eqref{2d_temporal} is similar to that in \cite{lauritzen2010conservative}.
We first present the algorithm for the second order $P^1$ SLDG scheme with a quadrilateral approximation for the shape of upstream cells. Then we generalize it to the third order $P^2$ SLDG scheme with quadratic-curved quadrilateral upstream cells, by highlighting the new components in the algorithm design.

\subsection{A $P^1$ SLDG Scheme with Quadrilateral Upstream Cells}
\label{sec:p1}

Below, we present the procedure of the proposed $P^1$ SLDG scheme with quadrilateral upstream cells. The algorithm consists of two main components: one is the search algorithm that finds overlapping regions between the upstream cell and background grid cells, i.e., Eulerian cells; and the other is the use of Green's theorem that enables us
 to convert the area integrals to line integrals when evaluating the RHS of \eqref{2d_temporal}. Such a procedure is similar to that in \cite{lauritzen2010conservative}.
%to evaluate the volume integral on the RHS of \eqref{2d_temporal}

 \begin{description}
   \item[(1)] \textbf{\emph{Characteristics tracing.}} Locate four vertices of upstream cell $A_{j}^\star$: $v_1^\star$, $v_2^\star$, $v_3^\star$ and $v_4^\star$ by tracking the characteristics backward to time $t^n$, i.e., solving the characteristics equations,
        \begin{equation}
        \begin{cases}
        \frac{d x(t) }{dt} = a(x(t) ,y(t) ,t ),\\
        \frac{ d y(t) }{dt} = b(x(t) ,y(t) ,t ), \\
        x(t^{n+1} ) = x(v_q),\\
        y(t^{n+1} ) = y(v_q)
        \end{cases}
        \end{equation}
       starting from the four vertices of $A_{j}$: $v_q$, $q=1,\cdots,4$ (see Figure \ref{schematic_2d} (a)).
      For the $P^1$ SLDG scheme it is sufficient to approximate the upstream cell by a quadrilateral to retain second order accuracy, see $A_j^\star$ in Figure \ref{schematic_2d} (b).

   \item[(2)] \textbf{\emph{Least squares approximation of test function $\psi(x, y, t^n)$.}}  Approximate $\psi(x,y,t^n)$  over the upstream cell $A_{j}^\star$.
              Specifically, based on the fact that the solution of the adjoint problem \eqref{2d_adjoint} stays unchanged along characteristics, we have
              \begin{equation*}
              \psi( x( v_q^\star ) ,y( v_q^\star ) , t^n ) = \Psi( x(v_q), y(v_q) ),\ \ q=1,2,3, \text{and} \ 4.
              \end{equation*}
              Thus, we can reconstruct a linear function $\psi^\star(x,y)$ approximating $\psi(x,y,t^n)$ on $A_{j}^\star$ by a least squares strategy and denote it as $\psi^\star(x,y)$.
              %=\bP^1\psi(x,y,t^n)$.
              % as such least square approximation.

   \item[(3)] %\textbf{\emph{Search algorithm.}}
   Denote $A_{j,l}^\star$ as a non-empty overlapping region between the upstream cell $A_j^\star$ and the background grid cell $A_l$, i.e., $A_{j,l}^\star = A_{j}^\star \cap A_l$, $A_{j,l}^\star \neq \emptyset,\ l\in \varepsilon_j^\star=\{ l| A_{j,l}^\star \neq \emptyset \},$ see Figure \ref{schematic_2d} (b). Then, we can  approximate the RHS of \eqref{2d_temporal} as follows
      \begin{equation}
             \iint_{A_j^\star} u(x,y,t^{n} )\psi(x,y,t^{n} ) dxdy
           \approx
          \sum_{l\in \varepsilon_j^\star}^{  } \iint_{A_{j,l}^\star } u(x,y,t^n)\psi^\star(x,y)dxdy.
       \label{temp1}
      \end{equation}

   \item[(4)] \textbf{\emph{Line integral evaluation.}} Note that the integrands on the RHS of \eqref{temp1} are piecewise quadratic polynomials.
       By introducing two auxiliary function $P(x,y)$ and $Q(x,y)$ such that
       \begin{equation*}
       -\frac{\partial P }{\partial y } + \frac{\partial Q}{\partial x }  =  u(x,y,t^n)\psi^\star(x,y),
       \end{equation*}
   the area integral $ \iint_{A_{j,l}^\star } u(x,y,t^n)\psi^\star(x,y)dxdy  $ can be converted into line integrals via Green's theorem, i.e.,
   \begin{equation}
      \iint_{A_{j,l}^\star } u(x,y,t^n)\psi^\star(x,y)dxdy = \oint_{\partial A_{j,l}^\star}  Pdx + Qdy,
      \label{Green}
   \end{equation}
   see Figure \ref{schematic_2d} (b). Note that the choices of $P$ and $Q$ are not unique, but the value of the line integrals is independent of the choices. In the implementation, we follow the same procedure in \cite{lauritzen2010conservative} when choosing $P$ and $Q$.

 \end{description}

\begin{figure}[h]
\centering

\subfigure[]{
\begin{tikzpicture}[scale = 1.1]
    %%%%%%%%%%%%%%%%%%%%%%%%%%%%%%%%%%%%%%%%%%%%%%%%%%%%%%%%%%mesh
    \draw[black,thin] (0,0.5) node[left] {} -- (4,0.5)
                                        node[right]{};
    \draw[black,thin] (0,2.) node[left] {$$} -- (4,2)
                                        node[right]{};
    \draw[black,thin] (0,3.5) node[left] {$$} -- (4,3.5)
                                        node[right]{};
                                        %%%%%%%%%%%
    \draw[black,thin] (0.5,0) node[left] {} -- (0.5,4)
                                        node[right]{};
    \draw[black,thin] (2,0) node[left] {$$} -- (2,4)
                                        node[right]{};
    \draw[black,thin] (3.5,0) node[left] {$$} -- (3.5,4)
                                        node[right]{};

    %%%%%%%%%%%%%%%%%%%%%%%%%%%%%%%%%%%%%%%%%%%%%%%%%%%
    %%%%%%%%%%%%%%%%%%%%%%%%%%%%%%%%%%%%%%%%%%%%%%%%%%% Lagrangian element
    \fill [red] (1.,1) circle (2pt) node[above right,black] {};
    \fill [red] (3,1) circle (2pt) node[] {};
    \fill [red] (1,2.5) circle (2pt) node[below right] {$A_j^\star$} node[above left] {};
    \fill [red] (2.5,2.5) circle (2pt) node[] {};
   \usetikzlibrary{shapes.geometric}
  \node[fill,star,star points=4, star point ratio=.2] at (2,1) {};
  \node[fill,star,star points=4, star point ratio=.2] at (2,2.5) {};
  \node[fill,star,star points=4, star point ratio=.2] at (1,2) {};
  \node[fill,star,star points=4, star point ratio=.2] at (2.65,2) {};

     \draw (0.5+0.01,2-0.01) node[fill=white,below right] {};
%*********************************************************************
\draw [red,thick] (1,1)node[right,scale=1.0]{}
        to[out=20,in=150] (2,0.7) node[] {};

        \draw [red,thick] (2,0.7)node[right,scale=1.0]{}
        to[out=330,in=240] (3,1) node[] {};
%********************************************************************
             \draw [red,thick] (1,2.5)node[right,scale=1.0]{}
        to[out=310,in=90] (1.1,2) node[] {};
        \draw [red,thick] (1.1,2)node[right,scale=1.0]{}
        to[out=270,in=80] (1,1) node[] {};
%*******************************************************************

        \draw [red,thick] (1,2.5)node[right,scale=1.0]{}
        to[out=10,in=180] (2.5,2.5) node[] {};

        \draw [red,thick] (3,1)node[right,scale=1.0]{}
        to[out=80,in=280] (2.5,2.5) node[] {};
    %%%%%%%%%%%%%%%%%%%%%%%%%%%%%%%%%%%%%%%%%%%%%%%%%%%%%%%%%%%%%
\draw [thick] (1,1)-- (3,1) node[] {};
 \draw [thick] (1,2.5) -- (1,1) node[] {};
        \draw [thick] (1,2.5)--(2.5,2.5) node[] {};
        \draw [thick] (3,1)--(2.5,2.5) node[] {};

\end{tikzpicture}

}
\subfigure[]{
\begin{tikzpicture}[scale = 1.1]
    %%%%%%%%%%%%%%%%%%%%%%%%%%%%%%%%%%%%%%%%%%%%%%%%%%%%%%%%%%mesh
    \draw[black,thin] (0,0.5) node[left] {} -- (4,0.5)
                                        node[right]{};
    \draw[black,thin] (0,2.) node[left] {$$} -- (4,2)
                                        node[right]{};
    \draw[black,thin] (0,3.5) node[left] {$$} -- (4,3.5)
                                        node[right]{};
                                        %%%%%%%%%%%
    \draw[black,thin] (0.5,0) node[left] {} -- (0.5,4)
                                        node[right]{};
    \draw[black,thin] (2,0) node[left] {$$} -- (2,4)
                                        node[right]{};
    \draw[black,thin] (3.5,0) node[left] {$$} -- (3.5,4)
                                        node[right]{};

    %%%%%%%%%%%%%%%%%%%%%%%%%%%%%%%%%%%%%%%%%%%%%%%%%%%
    %%%%%%%%%%%%%%%%%%%%%%%%%%%%%%%%%%%%%%%%%%%%%%%%%%% Lagrangian element
    \fill [red] (1.,1) circle (2pt) node[above right,black] {};
    \fill [red] (3,1) circle (2pt) node[] {};
    \fill [red] (1,2.5) circle (2pt) node[below right] {} node[above left] {};
    \fill [red] (2.5,2.5) circle (2pt) node[] {};

     \draw (0.5+0.01,2-0.01) node[fill=white,below right] {};
%*********************************************************************

    %%%%%%%%%%%%%%%%%%%%%%%%%%%%%%%%%%%%%%%%%%%%%%%%%%%%%%%%%%%%%
\draw [thick] (1,1)-- (3,1) node[] {};
 \draw [thick] (1,2.5) -- (1,1) node[] {};
        \draw [thick] (1,2.5)--(2.5,2.5) node[] {};
        \draw [thick] (3,1)--(2.5,2.5) node[] {};
   \usetikzlibrary{shapes.geometric}
  \node[fill,star,star points=4, star point ratio=.2] at (2,1) {};
  \node[fill,star,star points=4, star point ratio=.2] at (2,2.5) {};
  \node[fill,star,star points=4, star point ratio=.2] at (1,2) {};
  \node[fill,star,star points=4, star point ratio=.2] at (2.65,2) {};
  \draw [-latex] (1-0.2,2) -- node[right=3pt]{$\mathcal{L}_q$}(1-0.2,1) node[] {};
  \draw [-latex] (1-0.2,2.5) -- (1-0.2,2) node[] {};

  \draw [-latex] (2.5,2.5+0.2) -- (2,2.5+0.2);
  \draw [-latex] (2,2.5+0.2) -- (1,2.5+0.2);

  \draw [-latex] (1,1-0.2) -- (2,1-0.2) node[] {};
  \draw [-latex] (2,1-0.2) -- (3,1-0.2) node[] {};
  \draw [-latex] (3+0.2,1) -- (2.65+0.2,2) node[] {};
  \draw [-latex] (2.65+0.2,2) -- (2.5+0.2,2.5) node[] {};

  %\node[draw,star,star points=6, star point ratio=.6] at (2,1) {};

%\usetikzlibrary{shapes.geometric}
%   \node[diamond,draw] at (0,0) {};
%   \node[trapezium,draw] at (2,0) {};
%   \node[semicircle,draw] at (4,0) {};
%   \node[star,draw] at (6,0) {};
%   \node[isosceles triangle,draw] at (8,0) {};
%   \node[circular sector,draw] at (10,0) {};
%   \node[cylinder,draw] at (12,0) {};

\end{tikzpicture}

}
\subfigure[]{
\begin{tikzpicture}[scale = 1.1]
%

%\draw (2,1) node[fill=white,below right] {};
%\draw (2,2.5) node[fill=white,above=2pt] {$s_2$};
%\draw (0.8,1.85) node[fill=white,below] {$s_3$};
%\draw (2.6,2) node[fill=white,above right] {$s_4$};
%\draw (2,2) node[fill=white,below right] {$c_1$};
\node [below right,blue] at (2,1) {$s_1$};
\node [above,blue] at (2,2.5) {$s_2$};
\node [below,blue] at (0.9,2.) {$s_3$};
\node [above right,blue] at (2.6,2) {$s_4$};
\node [below right, blue] at (2,2) {$c_1$};
    %%%%%%%%%%%%%%%%%%%%%%%%%%%%%%%%%%%%%%%%%%%%%%%%%%%%%%%%%%mesh
    \draw[black,thin] (0,0.5) node[left] {} -- (4,0.5)
                                        node[right]{};
    \draw[black,thin] (0,2.) node[left] {$$} -- (4,2)
                                        node[right]{};
    \draw[black,thin] (0,3.5) node[left] {$$} -- (4,3.5)
                                        node[right]{};
                                        %%%%%%%%%%%
    \draw[black,thin] (0.5,0) node[left] {} -- (0.5,4)
                                        node[right]{};
    \draw[black,thin] (2,0) node[left] {$$} -- (2,4)
                                        node[right]{};
    \draw[black,thin] (3.5,0) node[left] {$$} -- (3.5,4)
                                        node[right]{};

    %%%%%%%%%%%%%%%%%%%%%%%%%%%%%%%%%%%%%%%%%%%%%%%%%%%
    %%%%%%%%%%%%%%%%%%%%%%%%%%%%%%%%%%%%%%%%%%%%%%%%%%% Lagrangian element
    \fill [red] (1.,1) circle (2pt) node[above right,black] {};
    \fill [red] (3,1) circle (2pt) node[] {};
    \fill [red] (1,2.5) circle (2pt) node[below right] {} node[above left] {};
    \fill [red] (2.5,2.5) circle (2pt) node[] {};

   \usetikzlibrary{shapes.geometric}
  \node[fill,star,star points=4, star point ratio=.2,blue] at (2,1) {};
  \node[fill,star,star points=4, star point ratio=.2,blue] at (2,2.5) {};
  \node[fill,star,star points=4, star point ratio=.2,blue] at (1,2) {};
  \node[fill,star,star points=4, star point ratio=.2,blue] at (2.65,2) {};
  \node[fill,star,star points=4, star point ratio=.2,blue] at (2,2) {};

  \draw [-latex] (2-0.1,2+0.1) -- (2-0.1,2.5-0.1) node[] {};
  \draw [-latex] (2+0.1,2.5-0.1)--(2+0.1,2+0.1 )  node[] {};

   \draw [-latex] (2-0.1,1+0.1) --node[auto]{$\mathcal{S}_q$} (2-0.1,2-0.1) node[] {};
  \draw [-latex] (2+0.1,2-0.1 )--(2+0.1,1+0.1 )  node[] {};
  \draw [-latex] (2.65-0.1,2-0.1) --(2+0.1,2-0.1)  node[] {};
  \draw [-latex] (2+0.1,2+0.1)--(2.65-0.1,2+0.1)  node[] {};
    \draw [-latex] (1+0.1,2+0.1) --(2-0.1,2+0.1)  node[] {};
  \draw [-latex] (2-0.1,2-0.1)--(1+0.1,2-0.1)  node[] {};

     \draw (0.5+0.01,2-0.01) node[fill=white,below right] {};
%*********************************************************************

    %%%%%%%%%%%%%%%%%%%%%%%%%%%%%%%%%%%%%%%%%%%%%%%%%%%%%%%%%%%%%
\draw [thick] (1,1)-- (3,1) node[] {};
 \draw [thick] (1,2.5) -- (1,1) node[] {};
        \draw [thick] (1,2.5)--(2.5,2.5) node[] {};
        \draw [thick] (3,1)--(2.5,2.5) node[] {};

\end{tikzpicture}

}
\caption{Schematic illustration of the search algorithm. $P^1$ case.}
\label{schematic_search}
\end{figure}
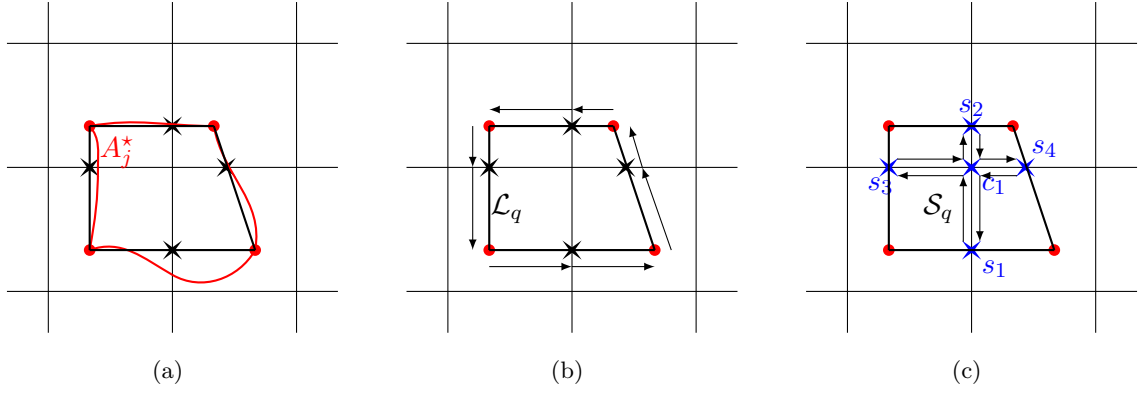
 In summary, combining \eqref{temp1} and \eqref{Green}, we have the following
\begin{align}
&
\iint_{A_{j}^\star } u(x,y,t_{n} )\psi(x,y,t_{n} ) dxdy  \notag \\
=&
\sum_{l\in \varepsilon_j^\star}^{  }\iint_{A_{j,l}^\star } u(x,y,t_{n} )\psi^\star(x,y ) dxdy   \notag \\
=&
\sum_{l\in \varepsilon_j^\star }^{  } \oint_{\partial A_{j,l}^\star}  Pdx + Qdy \notag \\
=&
\sum_{q=1}^{N_o}
\int_{ \mathcal{L}_q } [P  dx +Q  d y ]  + \sum_{q=1}^{N_i}
\int_{ \mathcal{S}_q } [Pdx +Q d y ] ,
\label{line}
\end{align}
where we compute the line integrals by organizing them as two categories: outer line segments (see Figure \ref{schematic_search} (b)) and inner line segments (see Figure \ref{schematic_search} (c)).
Again the procedure in evaluating the volume integral is the same as that in CSLAM and we refer to \cite{lauritzen2010conservative} for more details.

For implementation, the \textbf{\emph{search algorithm of line segments}} is provided as follows. We compute all intersection points of the four sides of the upstream cell with grid lines and organize them in the counterclockwise order. Then, we connect the organized points in the counterclockwise orientation and obtain the outer line segments, denoted as $\mathcal{L}_q$, $q=1,\cdots,N_o$, see Figure \ref{schematic_search} (b).
  The line-segments that are aligned with grid lines and enclosed by $A_j^\star$ are defined as super inner line segments.
  For example, in Figure \ref{schematic_search} (c), there are two super inner line segments, denoted by $s_1 s_2$ and $s_3 s_4$.
  We can thus find the inner line segments, denoted by $\mathcal{S}_q$, $q=1,\cdots,N_i$,  in each super line segment. In particular, four inner line segments $\overrightarrow{s_1 c_1}$, $\overrightarrow{c_1s_2}$, $\overrightarrow{s_2 c_1}$ and $\overrightarrow{c_1s_1}$ are obtained through breaking the super line segment $s_1 s_2$. Note that the orientation of inner segments must be taken in account when evaluating the associated line integrals. For instance, $\overrightarrow{s_1c_1}$ belongs to the left background cell and $\overrightarrow{c_1s_1}$ belongs to the right background cell.
  The indices of inner line segments in the background cell are ordered in a counterclockwise way.

 \begin{rem}
 Note that the above procedure can directly generalize to $P^k$ SLDG schemes with quadrilateral upstream cells, while the use of quadrilateral approximation will yield the second order accuracy for any high order $P^k$ approximation spaces. We will test the performance of such a configuration with $k=2$ in Section 4. Such observation motivates the construction of quadratic-curved quadrilateral approximation to the upstream cell in the next subsection for a third order $P^2$ SLDG scheme.
 \end{rem}

  \begin{rem}
 A similar second order scheme named characteristic discontinuous Galerkin method with quadrilateral upstream cells was proposed in \cite{lee2016high}. This method approximates a quadrilateral upstream cell by breaking quadrilateral upstream cells into several triangles and then generating quadrature points over each triangle. As mentioned, the generalization of  the strategy to the third order accuracy is very demanding in implementation. On the other hand,
  	our scheme can be naturally extended to the third order accuracy via the use of quadratic-curved quadrilateral upstream cells.
  \end{rem}

\label{section:p1}
 \subsection{A $P^2$ SLDG Scheme with Quadratic-curved Quadrilateral Upstream Cells}

We now present the  procedure of the $P^2$ SLDG scheme with quadratic-curved quadrilateral upstream cells to achieve a formal third order accuracy.
In particular, we propose to construct a parabola in approximating four sides of the upstream cell.  Since the scheme formulation is similar to that of the $P^1$ SLDG scheme described in the preceding subsection, we only highlight new ingredients in the algorithm design.

\begin{figure}[h]
\centering
\subfigure[]{
\begin{tikzpicture}
    %%%%%%%%%%%%%%%%%%%%%%%%%%%%%%%%%%%%%%%%%%%%%%%%%%%%%%%%%%mesh
    \draw[black,thin] (0,0.5) node[left] {} -- (5.5,0.5)
                                        node[right]{};
    \draw[black,thin] (0,2.) node[left] {$$} -- (5.5,2)
                                        node[right]{};
    \draw[black,thin] (0,3.5) node[left] {$$} -- (5.5,3.5)
                                        node[right]{};
    \draw[black,thin] (0,5 ) node[left] {$$} -- (5.5,5)
                                        node[right]{};
                                        %%%%%%%%%%%
    \draw[black,thin] (0.5,0) node[left] {} -- (0.5,5.5)
                                        node[right]{};
    \draw[black,thin] (2,0) node[left] {$$} -- (2,5.5)
                                        node[right]{};
    \draw[black,thin] (3.5,0) node[left] {$$} -- (3.5,5.5)
                                        node[right]{};
    \draw[black,thin] (5,0) node[left] {$$} -- (5,5.5)
                                        node[right]{};
    %%%%%%%%%%%%%%%%%%%%%%%%%%%%%%%%%%%%%%%%%%% background element
    \fill [blue] (3.5,3.5) circle (2pt) node[] {};
    \fill [blue] (5,3.5) circle (2pt) node[] {};
    \fill [blue] (3.5,5) circle (2pt) node[below right] {$A_j$} node[above left] {$v_7$};
    \fill [blue] (5,5) circle (2pt) node[] {};

    \fill [blue] (3.5,4.25) circle (2pt) node[] {};
    \fill [blue] (5,4.25) circle (2pt) node[] {};
    \fill [blue] (4.25,4.25) circle (2pt) node[] {};
        \fill [blue] (4.25, 3.5) circle (2pt) node[] {};
    \fill [blue] (4.25,5) circle (2pt) node[] {};

     \draw[thick,blue] (3.5,3.5) node[left] {} -- (3.5,5)
                                        node[right]{};
      \draw[thick,blue] (3.5,3.5) node[left] {} -- (5,3.5)
                                        node[right]{};
       \draw[thick,blue] (3.5,5) node[left] {} -- (5,5)
                                        node[right]{};
        \draw[thick,blue] (5,3.5) node[left] {} -- (5,5)
                                        node[right]{};
    %%%%%%%%%%%%%%%%%%%%%%%%%%%%%%%%%%%%%%%%%%%%%%%%%%%
    %%%%%%%%%%%%%%%%%%%%%%%%%%%%%%%%%%%%%%%%%%%%%%%%%%% Lagrangian element
    \fill [red] (1.,1) circle (2pt) node[above right,black] {$A_{j,l}^\star$};
    \fill [red] (3,1) circle (2pt) node[] {};
    \fill [red] (1,2.5) circle (2pt) node[below right] {$A_j^\star$} node[above left] {$v_7^\star$};
    \fill [red] (2.5,2.5) circle (2pt) node[] {};

     \draw[-latex,dashed](3.5,5)node[right,scale=1.0]{}
        to[out=240,in=70] (1,2.50) node[] {};

     \draw (0.5+0.01,2-0.01) node[fill=white,below right] {$A_l$};

     \draw [red,thick] (1,1)node[right,scale=1.0]{}
        to[out=20,in=150] (2,0.7) node[] {};

        \draw [red,thick] (2,0.7)node[right,scale=1.0]{}
        to[out=330,in=240] (3,1) node[] {};
%********************************************************************
             \draw [red,thick] (1,2.5)node[right,scale=1.0]{}
        to[out=310,in=90] (1.1,2) node[] {};
        \draw [red,thick] (1.1,2)node[right,scale=1.0]{}
        to[out=270,in=80] (1,1) node[] {};

        \draw [red,thick] (1,2.5)node[right,scale=1.0]{}
        to[out=10,in=180] (2.5,2.5) node[] {};

        \draw [red,thick] (3,1)node[right,scale=1.0]{}
        to[out=80,in=280] (2.5,2.5) node[] {};
    %%%%%%%%%%%%%%%%%%%%%%%%%%%%%%%%%%%%%%%%%%%%%%%%%%%%%%%%%%%%%

       \fill [red] (2,0.7) circle (2pt) node[above right,black] {};
    \fill [red] (2,1.7) circle (2pt) node[] {};
    \fill [red] (1.9,2.5) circle (2pt) node[below right] {} node[above left] {};
    \fill [red] (1.1,1.8) circle (2pt) node[] {};
\fill [red] (2.8,1.8) circle (2pt) node[] {};

\end{tikzpicture}
}
\subfigure[]{

\begin{tikzpicture}[scale = 1.3]
    %%%%%%%%%%%%%%%%%%%%%%%%%%%%%%%%%%%%%%%%%%%%%%%%%%%%%%%%%%mesh
    \draw[black,thin] (0,0.5) node[left] {} -- (4,0.5)
                                        node[right]{};
    \draw[black,thin] (0,2.) node[left] {$$} -- (4,2)
                                        node[right]{};
    \draw[black,thin] (0,3.5) node[left] {$$} -- (4,3.5)
                                        node[right]{};
                                        %%%%%%%%%%%
    \draw[black,thin] (0.5,0) node[left] {} -- (0.5,4)
                                        node[right]{};
    \draw[black,thin] (2,0) node[left] {$$} -- (2,4)
                                        node[right]{};
    \draw[black,thin] (3.5,0) node[left] {$$} -- (3.5,4)
                                        node[right]{};

    %%%%%%%%%%%%%%%%%%%%%%%%%%%%%%%%%%%%%%%%%%%%%%%%%%%
    %%%%%%%%%%%%%%%%%%%%%%%%%%%%%%%%%%%%%%%%%%%%%%%%%%% Lagrangian element
    \fill [red] (1.,1) circle (2pt) node[above right,black] {};
    \fill [red] (3,1) circle (2pt) node[] {};
    \fill [red] (1,2.5) circle (2pt) node[below right] {} node[above left] {};
    \fill [red] (2.5,2.5) circle (2pt) node[] {};

     \draw (0.5+0.01,2-0.01) node[fill=white,below right] {};

\draw [red,thick] (1,1)node[right,scale=1.0]{}
        to[out=20,in=150] (2,0.7) node[] {};

        \draw [red,thick] (2,0.7)node[right,scale=1.0]{}
        to[out=330,in=240] (3,1) node[] {};
%********************************************************************
             \draw [red,thick] (1,2.5)node[right,scale=1.0]{}
        to[out=310,in=90] (1.1,2) node[] {};
        \draw [red,thick] (1.1,2)node[right,scale=1.0]{}
        to[out=270,in=80] (1,1) node[] {};

        \draw [red,thick] (1,2.5)node[right,scale=1.0]{}
        to[out=10,in=180] (2.5,2.5) node[] {};

        \draw [red,thick] (3,1)node[right,scale=1.0]{}
        to[out=80,in=280] (2.5,2.5) node[] {};

               \fill [red] (2,0.7) circle (2pt) node[above right,black] {};
    \fill [red] (2,1.7) circle (2pt) node[] {};
    \fill [red] (1.9,2.5) circle (2pt) node[below right] {} node[above left] {};
    \fill [red] (1.1,1.8) circle (2pt) node[] {};
\fill [red] (2.8,1.8) circle (2pt) node[] {};

%%%%%%%%%%%%%%%%%%%%%%%%%%%%%%%%%%%%%%%%%%%%%%%%%%%%%%%%%%%%%%%coordinate
\draw[-latex,blue,thick](0,1)node[right,scale=1.0]{}
        to (4,1) node[below] {$\xi$};
\draw[-latex,blue,thick](2,0)node[right,scale=1.0]{}
        to (2,4) node[right] {$\eta$};

 %%%%%%%%%%%%%%%%%%%%%%%%%%%%%%%%%%%%%%%%%%%%%%%%%%%%%%%%%%%
       \draw[black, ultra thick] (2,0.7) node[below left =2pt] {$v_2^\star$}
        parabola(1,1)node[above left,scale=1.0]{ $v_1^\star$ };
   \draw[-latex,black, ultra thick](2,0.7)node[above left,scale=1.0]{  }
         parabola (3,1) node[below =2pt] {$v_3^\star$};
\end{tikzpicture}

}
\caption{Schematic illustration of the SLDG formulation in two dimensions. $P^2$ case.}
\label{schematic_2d_p2}
\end{figure}
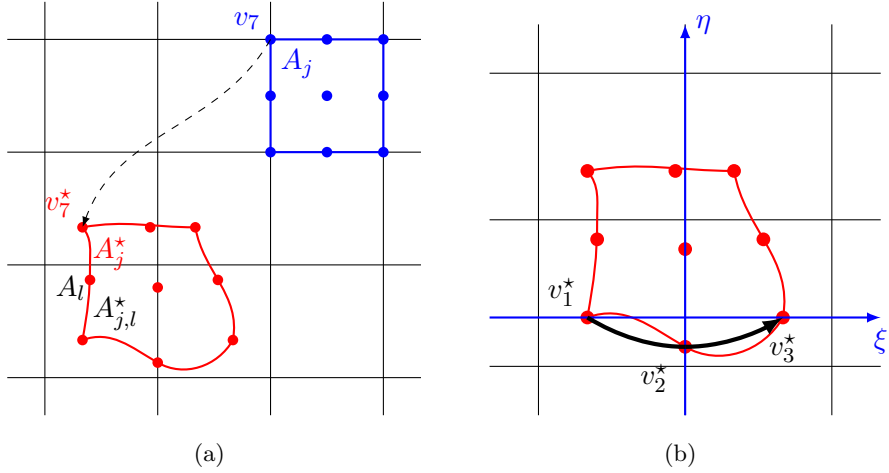

 \begin{description}
  \item[(1)] \textbf{\emph{Characteristics tracing.}} Locate the nine vertices of upstream cell $A_{j}^\star$:
  $v_i^\star,\,i=1,\ldots,9$ by  tracking the characteristics backward to time $t^n$, i.e., solving the characteristics equations, for the nine vertices of $A_{j}$:
  $v_i,\,i=1,\ldots,9$  (see Figure \ref{schematic_2d_p2} (a)).

   \item[(2)] \textbf{\emph{Reconstructing sides of upstream cells.}} Construct a quadratic curve to approximate each side of the upstream cell. In particular, to construct the quadratic curve, $\wideparen{ v_1^\star,v_2^\star, v_3^\star }$ as shown in Figure \ref{schematic_2d_p2} (b), we do the following.
   \begin{description}
     \item[(2a)]
     Construct  a coordinate transformation $x-y$ to $\xi-\eta$
      such that the coordinates of $v_1^\star$ and $v_3^\star$ are $(-1, 0)$ and $(1, 0)$ in $\xi-\eta$ space, respectively
      (see Figure \ref{schematic_2d_p2} (b)).
      Let $(x_1, y_1)$ and $(x_3, y_3)$ be the $x-y$ coordinate of $v_1^\star$ and $v_3^\star$, then
      the coordinate transformation is obtained as
     \begin{equation}
     \begin{cases}
     \xi(x,y) = a x + by + c,\\
     \eta(x,y) = bx -ay + d,
     \end{cases}
     \label{trans}
     \end{equation}
     where
     $$a=\frac{ 2( x_3 -x_1  ) }{ (x_1 - x_3  )^2 + (y_1 - y_3  )^2  },$$
     $$b=\frac{ 2( y_3 -y_1  ) }{ (x_1 - x_3  )^2 + (y_1 - y_3  )^2  },$$
     $$c=\frac{   x_1^2 - x_3^2 + y_1^2 - y_3^2 }{ (x_1 - x_3  )^2 + (y_1 - y_3  )^2  },$$
     $$d=\frac{ 2(x_3y_1 - x_1y_3  ) }{ (x_1 - x_3  )^2 + (y_1 - y_3  )^2  }.$$
     Its reverse transformation can be constructed accordingly:
     \begin{equation}
     \begin{cases}
     x = \frac{x_3 -x_1 }{2} \xi  + \frac{ y_3 -y_1}{2} \eta + \frac{x_3 +x_1}{2}, \\
     y = \frac{ y_3 -y_1}{2} \xi  - \frac{ x_3 -x_1 }{2} \eta + \frac{ y_3 + y_1}{2}.
     \end{cases}
     \label{reverse}
     \end{equation}

     \item[(2b)] Get the $\xi-\eta$ coordinate for the point $v_2^\star$ as $(\xi_2, \eta_2)$. Based on $(-1,0)$, $(\xi_2,\eta_2)$ and $(1,0)$, we construct the parabola,
         \begin{equation}
         \wideparen{ v_1^\star,v_2^\star, v_3^\star }: \eta(x,y) = \frac{\eta_2}{ \xi_2^2-1 }(  \xi(x,y)^2-1  ).
         \label{parabola}
         \end{equation}

   \end{description}

   \item[(3)] \textbf{\emph{Search algorithm of line segments.}} The scheme is implemented by  evaluating the line integrals along outer line segments and inner line segments in the RHS of \eqref{line}, similar to the quadrilateral case. To organize these outer and inner line segments, we perform the following procedure.

   \begin{description}
   \item[(3a)] Outer segments: compute all the intersections between grid lines and four curved-sides  of the upstream cell $A_{j}^\star$ and organize them in the  counterclockwise order to obtain the outer line segments $\mathcal{L}_q$, $q=1,\cdots,N_o$. To find intersection points, we solve the following equations

     \begin{equation}
     \begin{cases}
     x_i =   \frac{x_3 -x_1 }{2} \xi + \frac{y_3 -y_1 }{2} \eta
       + \frac{ x_3 +x_1 }{2} \   (\text{intersection with grid line}\  x=x_i ), \\
      \eta = \frac{\eta_2}{ \xi_2^2-1 }(  \xi^2-1  )
     \end{cases}
     \end{equation}
     and
          \begin{equation}
     \begin{cases}
     y_j  =   \frac{y_3 -y_1 }{2} \xi - \frac{x_3 -x_1 }{2} \eta
       + \frac{ y_3 +y_1 }{2} \   (\text{intersection with grid line}\  y=y_j ) ,\\
      \eta = \frac{\eta_2}{ \xi_2^2-1 }(  \xi^2-1  ) .
     \end{cases}
     \end{equation}
     If $\xi \in [-1,1]$, the solution $(\xi, \eta)$ is identified as an intersection point.

      \item[(3b)]   Inner segments: collect all intersections from four sides and follow a similar procedure provided in Section \ref{section:p1}, then the inner line segment $\mathcal{S}_q$, $q=1,\cdots,N_i$ can be determined.

   \end{description}

   \item[(4)] \textbf{\emph{Line integral evaluation.}}
    The integral of inner line segments $ \sum_{q=1}^{N_i}
\int_{ \mathcal{S}_q } [Pdx +Q d y ]$ can be evaluated in the same way as in the $P^1$ case.
The integral of outer line segments $\sum_{q=1}^{N_o}
\int_{ \mathcal{L}_q } [P  dx +Q  d y ]  $ can be evaluated by the following parameterization on each side.
Assume that $\mathcal{L}_q$ is the part of the side $\wideparen{v_1^\star v_2^\star v_3^\star}$.
Substitute \eqref{parabola}  into \eqref{reverse}, we have
     \begin{equation}
     \begin{cases}
     x(\xi) = \frac{x_3 -x_1 }{2} \xi  + \frac{ y_3 -y_1}{2} \frac{\eta_2}{ \xi_2^2-1 }(  \xi^2-1  ) + \frac{x_3 +x_1}{2}, \\
     y(\xi) = \frac{ y_3 -y_1}{2} \xi  - \frac{ x_3 -x_1 }{2} \frac{\eta_2}{ \xi_2^2-1 }(  \xi^2-1  ) + \frac{ y_3 + y_1}{2}.
     \end{cases}
     \end{equation}

Hence,
\begin{equation}
\int_{  \mathcal{L}_q } [Pdx+ Qdy ]
= \int_{\xi^{(q)} }^{\xi^{(q+1)} }  [P( x(\xi,\eta),y(\xi,\eta) ) x'(\xi)+ Q( (\xi,\eta),y(\xi,\eta) ) y'(\xi) ]d\xi,
\end{equation}
where $(\xi^{(q)},\eta^{(q)})$ and $(\xi^{(q+1)},\eta^{(q+1)})$  are the starting point and the end point of $\mathcal{L}_q$ in $\xi-\eta$ coordinate, respectively.
The above integral can be done analytically or by numerical quadrature rules.

 \end{description}

\begin{rem}
Note that the proposed SLDG schemes are mass conservative, as with CSLAM \cite{lauritzen2010conservative,erath2013mass}, if the boundary condition of \eqref{2d_linear} is periodic or has compact support.
Letting the test function $\psi=1$, the total mass on all upstream cells can be written as follows,
\begin{align}
\sum_{j=1}^J \iint_{A_j^\star} u(x,y,t^n) dxdy &= \sum_{j=1}^J \sum_{ l\in \varepsilon_j^\star } \iint_{A_{j,l}^\star} u(x,y,t^n)dxdy,
\end{align}
where the first summation is with respect to the index for upstream cells and the second summation is with respect to their nonempty overlapping regions with background Eulerian cells. On the other hand, denote $A_{j,l}$ as a non-empty overlapping region between the Eulerian cell $A_j$ and the upstream cell $A^*_l$.
The total mass on all upstream cells equals the total mass on all Eulerian cells under the periodic boundary condition, i.e.
\begin{align}
\label{eq: mass}
\sum_{j=1}^J \sum_{ l\in \varepsilon_j^\star } \iint_{A_{j,l}^\star} u(x,y,t^n)dxdy
=
\sum_{j=1}^J \sum_{ l\in \varepsilon_j } \iint_{A_{j,l} } u(x,y,t^n)dxdy
=
\sum_{j=1}^J \iint_{A_j  } u(x,y,t^n) dxdy,
\end{align}
where $\varepsilon_j= \{ l | A_{j,l}  \neq \emptyset   \}$. In the second expression in eq.~\eqref{eq: mass}, the first summation sign is with respect to the index for Eulerian cells while the second one is with respect to their nonempty overlapping regions with upstream cells.
Note that the area integral $\iint_{A_{j,l} } u(x,y,t^n)dxdy $ is converted to  line integrals in  the implementation of the proposed SLDG scheme.
The SLDG scheme maintains the mass conservation, given that an upstream cell shares the same side with its adjacent neighbors.
%This SLDG scheme is still  mass conservative, since
%the summation of all line integrals of outer segments (except when upstream cell sides coincide with grid lines) in a Eulerian cell $A_{j}$ yields zero.
%In fact, a
%line integral along a particular side of a upstream cell is exactly equal to the line integral along the same side shared with the adjacent upstream cell but with opposite sign.

\end{rem}

\begin{rem}
	To the best of authors' knowledge, the proposed methodology is the first non-splitting scheme that is able to attain the formal third order accuracy and allows for large time step evolution and mass conservation.
%Note that the above third order accuracy, where a quadratic-curved quadrilateral approximation is used, is the first attempt for high order approximation of upstream cell under the semi-Lagrangian framework as far as we know.
\end{rem}

\begin{rem}
	When the velocity field $\ba=(a,b)$ is constant, an upstream cell will have the same shape as the Eulerian cell. Based on the observation,  using the strategy in \cite{qiu2011positivity}, we can easily show that the proposed SLDG schemes are $L^2$ stable and establish a prior error estimate accordingly.  However, when   $\ba$ becomes space and time dependent, one needs to use either a quadrilateral or a quadratic-curved quadrilateral to approximate each upstream cell, see \eqref{temp1}. As such, the Galerkin weak formulation \eqref{2d_temporal} is not exactly computed but subject to some approximation error,  leading to a variational crime (see \cite{brenner2008mathematical}, Chapter 10). In our future work, we will study such a variational crime and establish the $L^2$ stability analysis in the general setting.
\end{rem}

\begin{rem}
	The proposed method was established based on a Cartesian partition of a 2D rectangular domain for simplicity. In principle, the method can be extended to unstructured meshes. The only modification needed lies in the search algorithm. In particular, instead of considering the intersection of each upstream cell with background Cartesian grid lines, we need to consider the intersection of each upstream cell (approximated by triangles or quadratic-curved triangles) and the sides of background Eulerian triangle cells, which may be more involved in the implementation. The extension of the algorithm to three dimensions or higher, by using the Green's Theorem, is subject to future investigation.
\end{rem}

\subsection{Bound-preserving (BP) Filter}

 With the assumption that the velocity field $\ba$ is divergence free,
if the initial condition for \eqref{eq:trans} is positive, then the solution always stays positive as time evolves. Such a property is called  positivity preserving. In our SLDG schemes, it can be shown that the updated cell averages at $t^{n+1}$ stay positive, if the cell averages at $t^n$ are positive.
Similar to \cite{qiu2011positivity,Guo2013discontinuous,guo2015efficient}, in order to preserve positivity of numerical solutions, we further apply a high order BP filter \cite{zhangshu2010} into the proposed SLDG
scheme, which can be implemented as follows. The numerical solution $u(x,y,t^n)$ in cell $A_j$ is modified by $\widetilde{u}(x,y)$
\begin{equation*}
\widetilde{u}(x,y) = \theta ( u(x,y,t^n)- \overline{u} ) + \overline{u},\
 \theta = \min \left\{  \left| \frac{\overline{u}}{m'-\overline{u}} \right|,1  \right\},
\end{equation*}
where $\overline{u}$ is the cell average of the numerical solution and $m'$ is the minimum value of $u(x,y,t^n)$ over $A_j$.
%In our 2D numerical tests, the minimum value of $P^k$ $(k=1,2)$ polynomial can be easily found.
 For $P^1$ polynomials, the minimum value can be found by comparing the values at four vertices of $A_j$.
 For $P^2$ polynomials, besides the four vertices, all critical points inside $A_j$ should be considered  to determine the function's minimum value. Note that the proposed SLDG schemes with the BP filter feature the $L^1$ conservation property hence the $L^1$ stability for nonnegative initial conditions, and the proof follows a similar argument in \cite{qiu2011positivity}.

\subsection{Data structure and flowchart}

In this subsection, we present the data structure and the flowchart of the proposed 2D SLDG schemes. The implementation is based on the object-oriented technology and the class hierarchy is shown in Figure \ref{class}. Compared with a structured approach, the object-oriented technology allows for  data encapsulation and hence greatly facilitates data access and manipulation.  Note that
a line in Figure \ref{class} represents a certain relationship between classes. The two main classes are \textbf{Cell-E} and \textbf{Cell-U} which represent  an Eulerian cell and an   upstream cell, respectively. There is a one-to-one relation between class \textbf{Cell-E} and class \textbf{Cell-U}. Recall that an upstream cell is obtained by tracing the characteristics backward for one time step. Class \textbf{Cell-E} has class \textbf{Node-E} and class \textbf{DG solution} as it data members: \textbf{Node-E} represents the vertices of the Eulerian cell and  class \textbf{DG solution} stores the DG solutions (coefficients).
Furthermore, class \textbf{Cell-U} has three data members including class \textbf{Node-U}, class \textbf{Test function}, and class \textbf{Segment}. In particular,
 class \textbf{Node-U} represents the vertices of the underlying upstream cell; class \textbf{Test function} stores the reconstructed test functions defined over the upstream cell; class \textbf{Segment} stores the data of the inner and outer segments. Also note that class \textbf{Side-U} that represents the sides of an upstream cell use class \textbf{Node-U} to define the start and end points. In addition, the intersection points of the sides of an upstream cell  (class \textbf{Side-U}) and the grid lines are stored in class \textbf{Intersection point}, which is a data member of class \textbf{Segment} and is used to obtain the inner and outer segments.

%  the  and \textbf{Side-U} represent the nodes and sides of the upstream cell.
%\textbf{Intersection point} class represents intersection points between the side of upstream cell and the side of Eulerian cell. In the SLDG scheme, test functions at $t^n$ , named \textbf{test function} class, are the part of the integrand in the area integral \eqref{temp1}. In fact, the area integral is converted to line integral evaluation of a lot segments including inner segments and outer segments on the upstream cell, formed \textbf{Segment} class.

Figure \ref{flowchart} presents the algorithm flowchart of the proposed SLDG schemes based on the class hierarchy introduced in Figure \ref{class}. In particular, such a flowchart summarizes the main procedures of the schemes developed in Section 3.1 and 3.2.
%After initialization, the DG solution on each Eulerian cell can be obtained.
% Firstly, adopting \emph{characteristics tracing} from \textbf{Node-E} class,  \textbf{Node-U} class can be obtained.
%Then,  a \emph{search algorithm} is used to obtain \textbf{Segment} class.
% Since the test function is a constant along the characteristic line, its value at a node can be obtained from the corresponding node on the Eulerian cell, and   the\textbf{ test function} $\psi^\star(x,y)$ can be reconstructed by \emph{least square approximation}.
% The area integral is converted to the line integral, and such line integral can be solved analytically.
% Finally, DG solution at $t^{n+1}$ can be updated via the summation of line integrals of \textbf{Segment} class.
%
%%%%%%%%%%%%%%%%%%%%%%%%%%%%%%%%%%%%%%%%%%%%%%%%%%%%%%%%%%%%%%%%%%%%%%%%%%%%%%%%%%%%%%%%%%%%%%%%%%%%%%%%%%%%%%
\tikzstyle{abstract}=[rectangle, draw=black,text centered, anchor=north, text width=2.5cm]
\tikzstyle{abstract1}=[rectangle, draw=black,anchor=north, text width=2.5cm]
\tikzstyle{comment}=[rectangle, draw=black, rounded corners, fill=green, drop shadow,
        text centered, anchor=north, text=white, text width=3cm]
\tikzstyle{myarrow}=[->, >=open triangle 90, thick]
\tikzstyle{compoarrow}=[->,>=diamond,thick ]
\tikzstyle{inhearrow}=[->,>=diamond,thick ]
\tikzstyle{line}=[-, thick]

\begin{figure}[h]
\begin{tikzpicture}[scale=0.9]
\tikzstyle{every node}=[font=\scriptsize]
    \node (EulerianCell) [abstract, rectangle split, rectangle split parts=1]
    at(0,2)
        {
            \textbf{Cell-E}
            %\nodepart{second}name
        };

      \node (Solution) [abstract, rectangle split, rectangle split parts=1]
    at(-2,0)
        {
            \textbf{DG solution}
            %\nodepart{second}name
        };
     \node (NodeE) [abstract, rectangle split, rectangle split parts=1]
    at(2,0)
        {
            \textbf{Node-E}
            %\nodepart{second}name
        };
   \draw[line] (Solution.north) -- ++(0,0.5) -| (EulerianCell.south);
     \draw[line] (NodeE.north) -- ++(0,0.5) -| (EulerianCell.south);

% \node (SideE) [abstract, rectangle split, rectangle split parts=1]
%    at(-2,-2)
%        {
%            \textbf{Side-E}
%            %\nodepart{second}name
%        };

   %\draw[line] (SideE.north) --  (NodeE.south);
   % \node (AuxNode01) [text width=3cm, right=of EulerianCell] {};
   %  \node (UpstreamCell)[abstract, rectangle split, rectangle split parts=1, right=of AuxNode01]
   \node (UpstreamCell)[abstract, rectangle split, rectangle split parts=1]
   at(10,2)
        {
            \textbf{Cell-U}
            %\nodepart{second}name
        };

     \draw[line] (EulerianCell.east) -- ++(0,0.0) -| (UpstreamCell.west);
%%%%%%%%%%%%%%%%%%%%%%%%%%%%%%%%%%%%%%%%%%%%%%%%%%%%%%%%%%%%%%%%%%%%%%%%%%%%%%%%%%%%%%%%%%%%%%%%%%%%%%%%%%%%%%%%

%\node (AuxNode02) [text width=2cm, below=of UpstreamCell] {};
% \node (Segment)[abstract, rectangle split, rectangle split parts=1, below left=of AuxNode02]
     \node (NodeU)[abstract, rectangle split, rectangle split parts=1]
     at(6,0)
        {
            \textbf{Node-U}
            %\nodepart{second}name
        };

  \node (SideU)[abstract, rectangle split, rectangle split parts=1]
     at(6,-2)
        {
            \textbf{Side-U}
            %\nodepart{second}name
        };

   % \node (NodeU)[abstract, rectangle split, rectangle split parts=1, below right=of AuxNode02]
   \node (Segment)[abstract, rectangle split, rectangle split parts=1]
   at(14,0)
        {
            \textbf{Segment}
            %\nodepart{second}name
        };
 \node (Crossing)[abstract, rectangle split, rectangle split parts=1]
   at(14,-1.75)
        {
            \textbf{Intersection point}
            %\nodepart{second}name
        };
     \draw[line] (NodeU.north) -- ++(0,0.5) -| (UpstreamCell.south);
     \draw[line] (Segment.north) -- ++(0,0.5) -| (UpstreamCell.south);
 %\node (AuxNode03) [text width=0.01pt, right=of Segment] {};
%  \node (Test)[abstract, rectangle split, rectangle split parts=1, right =of Segment]
  \node (Test)[abstract, rectangle split, rectangle split parts=1]
  at(10,0)
        {
            \textbf{Test function}
            %\nodepart{second}name
        };
     \draw[line] (Test.north) --  (UpstreamCell.south);

     \draw[line] (Crossing.north) --  (Segment.south);
     \draw[line] (SideU.north) --  (NodeU.south);

     %\draw[line] (SideE.east) --  (Crossing.west);
     \draw[line] (SideU.east) --  (Crossing.west);

     \draw[line] (NodeE.east) --  (NodeU.west);

     %%%%%%%%%%%%%%%%%%%%%%%%%%%%%%%%%%%%%%%%%%%%%%%%%%%%%%%%
     %%%%%%%%%%%%%%%%%%%%%%%%%%%%%%%%%%%%%%%%%%%%%%%%%%%%%%%%

\fill [black] (-3,-3.5-0.5) circle (2pt) node[right=2pt] {Cell-E: Eulerian Cell.};
\fill [black] (-3,-4.-0.5) circle (0pt) node[right=2pt] {Node-E: $\{v_q\}_{1}^4$ for $P^1$ or $\{v_q\}_{1}^9$ for $P^2$.  };
%\fill [black] (-3,-4.5-0.5) circle (0pt) node[right=2pt] {Side-E: Side of Eulerian Cell. };

\fill [black] (6,-3.5-0.5) circle (2pt) node[right=2pt] {Cell-U: Upstream Cell.};
\fill [black] (6,-4-0.5) circle (0pt) node[right=2pt] {Node-U: $\{v_q^\star\}_{1}^4$ for $P^1$ or $\{v_q^\star\}_{1}^9$ for $P^2$.  };
\fill [black] (6,-4.5-0.5) circle (0pt) node[right=2pt] {Side-U: Side of Upstream Cell.};

\fill [black] (-3,-5 ) circle(0pt) node[right=2pt] {DG solution: $u(x,y,t^n)$ . };
\fill [black] (-3,-5.5-0.5) circle(2pt) node[right=2pt] {Intersection point: Intersecion points of the sides of an upstream cell and grid lines. };

\fill [black] (6,-5-0.5) circle(0pt) node[right=2pt] {Test function: $\psi^\star(x,y)$ . };

\end{tikzpicture}
\caption{Class diagram for SLDG code.}
\label{class}
\end{figure}

\tikzstyle{myarrow}=[->, >=open triangle 90, thick]
\tikzstyle{compoarrow}=[->,>=diamond,thick ]
\tikzstyle{inhearrow}=[->,>=diamond,thick ]
\tikzstyle{line}=[-, thick]
\tikzstyle{arrowline} = [draw, -latex']
\tikzstyle{arrowline1} = [draw, latex'-]
\tikzstyle{decision} = [diamond, draw=black, text width=1.5cm]

\begin{figure}[h]
\centering
\begin{tikzpicture}[scale=0.9]

\tikzstyle{every node}=[font=\scriptsize]
    \node (EulerianCell) [abstract1, rectangle split, rectangle split parts=1]
    at(-4,0.55+0.5)
        {
        	\textbf{Initialization}:
            set parameters, \\
            generate meshes, \\
            and initialize the solution at $t^0$.

            %\nodepart{second}name
        };

    \node (DGsolution) [abstract1, rectangle split, rectangle split parts=1]
    at(0,0+0.5)
        {
            \textbf{
            DG solution} at $t^n$

            %\nodepart{second}name
        };

   \node (Chara) [abstract1, rectangle split, rectangle split parts=1]
    at(5,2)
        {

            Compute \textbf{Node-U} from \textbf{Node-E}
            via \textbf{\emph{characteristic tracing}}.

            %\nodepart{second}name
        };

       \node (SideU) [abstract1, rectangle split, rectangle split parts=1]
    at(3,-1)
        {

            Compute \textbf{Side-U} from \textbf{Node-U}

            %\nodepart{second}name
        };
       \node (Test) [abstract1, rectangle split, rectangle split parts=1]
    at(7,-1)
        {

            Approximate \textbf{test function} $ \psi^\star(x,y) $ via \textbf{\emph{least squares}}.

            %\nodepart{second}name
        };

   \node (Segment) [abstract1, rectangle split, rectangle split parts=1]
    at(3,-2.5)
        {

            Compute \textbf {intersection point} and then \textbf{Segment} via the \textbf{\emph{search algorithm}}.

            %\nodepart{second}name
        };

 \node (Line) [abstract1, rectangle split, rectangle split parts=1]
    at(5,-6)
        {
        	Obtain
            \textbf{
            \emph{Line integrals} via \textbf{\emph{Green's theorem}}
            }
            %\nodepart{second}name
        };

   \node (Update) [abstract1, rectangle split, rectangle split parts=1]
    at(5,-7.5)
        {

            Update \textbf{DG solution
            }.
            %\nodepart{second}name
        };

     \node (Decision) [decision] at(0,-8)
     {$t^{n+1}< Tfinal$ };

     \fill [black] (-1.8,-1.5) circle(0pt) node[right=0pt] {$n=n+1$};

\fill [black] (-0.8,-6.3) circle(0pt) node[right=2pt] {Yes};
  \fill [black] (-2.1,-7.6) circle(0pt) node[right=2pt] {No};
        \node (Stop) [abstract, rectangle split, rectangle split parts=1]
    at(-4,-7.6)
        {

            Stop

            %\nodepart{second}name
        };

\draw[arrowline](DGsolution.north) -- ++(0,0.25) |- (Chara.west) ;

       \draw[arrowline](Decision.west) --(Stop.east) ;

      \draw[arrowline](Update.west) --(Decision.east);

    %\draw[arrowline](EulerianCell.south) --(DGsolution.north);
    \draw[arrowline](EulerianCell.east) --(DGsolution.west);

    \draw[arrowline](Decision.north) -- (DGsolution.south) ;

     %\draw[arrowline](Segment.south) -- ++(0,-0.2) -|(Line.north);
     \draw[arrowline1](Line.north) -- ++(0,0.3) -|(Segment.south);
     \draw[arrowline1](Line.north) -- ++(0,0.3) -|(Test.south);
     \draw[arrowline1](Line.north) -- ++(0,0.3) -|(0.3,-0.9+0.5);

      \draw[arrowline](Line.south) --(Update.north);
       \draw[arrowline](SideU.south) --(Segment.north);

    \draw[arrowline](Chara.south) -- ++(0,-0.3) -|(SideU.north);
    \draw[arrowline](Chara.south) -- ++(0,-0.3) -|(Test.north);
\end{tikzpicture}
\caption{Flowchart for SLDG code.}
\label{flowchart}
\end{figure}
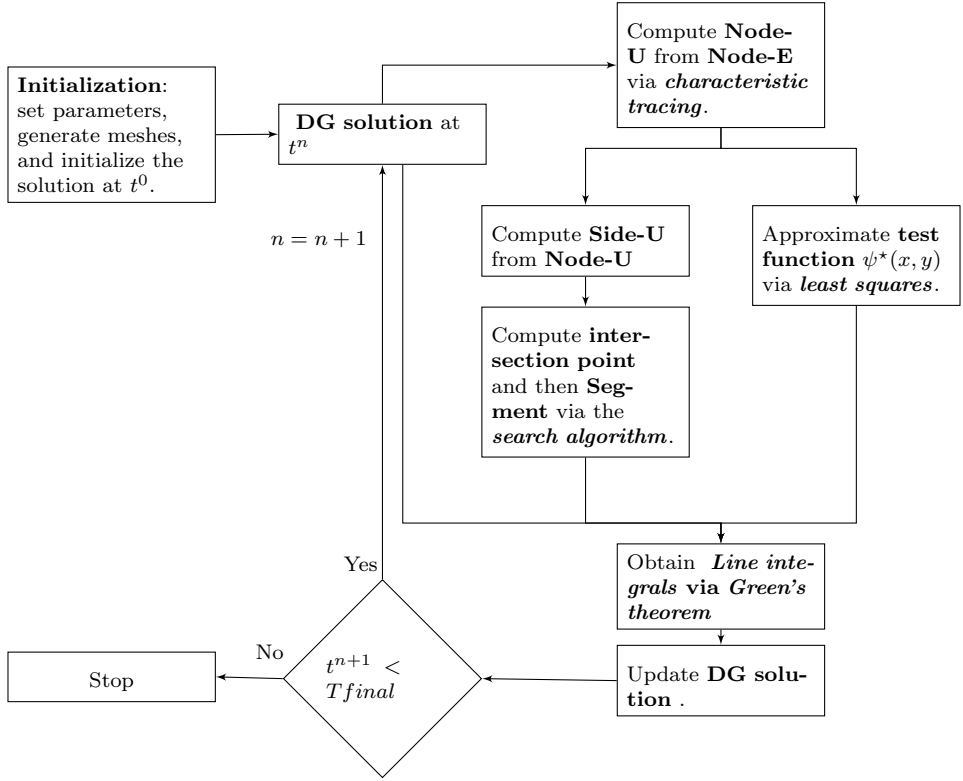

%\subsection{A Weak Formulation}
%
%\subsection{Algorithm Flowchart}
%
%\subsubsection{Search Algorithm}
%
%\subsubsection{Converting Area Integrals into Line Integrals}

%\input{analysis}
\section{Numerical Results}

In this section, we first test the performance of the $P^k$ SLDG schemes ($k=1,2,3$) for 1D examples. Then, we perform numerical tests in 2D for the $P^k$ SLDG schemes ($k=1,2$) with quadrilateral upstream cell (denoted by $P^k$ SLDG) and $P^2$ SLDG with quadratic-curved quadrilateral upstream cells (denoted by $P^2$ SLDG-QC).
In particular, for the constant-coefficient transport and the rigid body rotation problems, we only test the $P^k$ SLDG scheme, since the sides of the upstream cells are always straight lines. For the swirling flow problem, we run numerical tests for both the $P^k$ SLDG with $k=1,\,2$ and $P^2$ SLDG-QC schemes. Note that the upstream cells are no longer quadrilateral and hence different approximations used will lead to different numerical performance.
Furthermore, similar to other DG methods, the proposed SLDG schemes will generate spurious oscillations if the underlying solutions involve discontinuities due to the famous Gibbs phenomenon.
In the simulations, we use a simple WENO limiter proposed in \cite{zhong2013simple} to suppress such undesired oscillations. The positivity-preserving property of the SLDG solutions can be realized by incorporating the BP limiter.

\begin{exa}
(1D linear equation.) We start with the following 1D transport equation
\begin{equation}
u_t + u_x = 0,\ x\in[0,2\pi],
\label{linear}
\end{equation}
with smooth initial data $u(x,0)=\sin(x)$ and exact solution $u(x,t)=\sin(x-t)$.
Table \ref{table:1dlinear} summarizes the convergence study including the $L^2$ and $L^\infty$ errors and the associated orders of accuracy of numerical solutions by $P^k$ SLDG methods with $\Delta t = 0.5\Delta x$ and $\Delta t = 2.5\Delta x$. The expected orders of convergence are observed. The errors from $\Delta t = 2.5\Delta x$ are either comparable or slightly smaller than those from simulations with $\Delta t = 0.5\Delta x$.

\begin{table}[!ht]%\small
\caption{SLDG for \eqref{linear} with $u(x,0)=\sin(x)$. $T=20$.}
\vspace{0.1in}
\centering
\begin{tabular}{l ll ll ll ll}

\hline
 mesh & \multicolumn{4}{c}{$P^k$ SLDG with $\Delta t= 0.5\Delta x $}
     &\multicolumn{4}{c}{$P^k$ SLDG with $\Delta t= 2.5\Delta x $ } \\
   \cmidrule(lr){2-5} \cmidrule(lr){6-9}
{   }  &{$L^2$ error} & Order  & {$L^{\infty}$ error} & Order  &{$L^2$ error} & Order  & {$L^{\infty}$ error} & Order \\
\hline
    \multicolumn{9}{l}{ $P^1$ SLDG}
     \\

    20 &     3.89E-03 &          &    7.66E-03 &
     & 3.20E-03 &   & 1.00E-02 &  \\
    40 &     8.08E-04 &     2.27 &     2.13E-03 &     1.85
     &  7.59E-04 &     2.08 &     2.46E-03 &     2.03 \\
    80 &     2.00E-04 &     2.01 &     6.18E-04 &     1.79
     &     1.97E-04 &     1.94 &     6.58E-04 &     1.90 \\
    160 &     4.89E-05 &     2.04 &     1.58E-04 &     1.97
     &     5.07E-05 &     1.96  &     1.71E-04 &     1.94 \\
     320 &     1.15E-05 &     2.09 &    3.73E-05 &     2.09
    &     1.15E-05 &     2.14 &     3.79E-05 &     2.18  \\
    \multicolumn{9}{l}{ $P^2$ SLDG}
     \\

    20 &     7.37E-05 &          &    1.42E-04 &
     & 7.37E-05 &   &1.42E-04 &  \\
    40 &     9.22E-06 &     3.00 &     2.98E-05 &     2.25
     &9.22E-06 &     3.00 &     2.98E-05 &     2.26\\
    80 &    1.15E-06 &     3.01 &     2.08E-06 &     3.84
     &     1.15E-06 &     3.01 &     2.08E-06 &     3.84\\
    160 &     1.44E-07 &     2.99 &     4.42E-07 &     2.24
     &      1.44E-07 &     2.99 &     3.27E-07 &     2.67  \\
     320 &     1.74E-08 &     3.05 &     2.98E-08 &     3.89
    &     1.74E-08 &     3.05 &     2.99E-08 &     3.45  \\

    \multicolumn{9}{l}{ $P^3$ SLDG }
     \\

    20 &     1.46E-06 &          &    3.30E-06 &
     & 1.46E-06 &   &3.32E-06 &  \\
    40 &     9.08E-08 &     4.00 &     2.38E-07 &     3.79
     &9.08E-08 &     4.00 &     2.38E-07 &     3.80\\
    80 &     5.80E-09 &     3.97 &     1.46E-08 &     4.03
     &     5.80E-09 &     3.97 &     1.46E-08 &     4.03\\
    160 &     3.52E-10 &     4.04 &     8.31E-10 &     4.13
     &     3.42E-10 &     4.08 &    6.48E-10 &     4.49 \\
     320 &     2.47E-11 &     3.83 &     8.12E-11 &     3.36
    &     2.47E-11 &     3.79 &     8.11E-11 &     3.00  \\  \hline
\end{tabular}
\label{table:1dlinear}
\end{table}

\end{exa}

\begin{exa}
(1D transport equation with variable coefficients.) Consider
\begin{equation}
u_t + (\sin(x) u )_x = 0,\ x\in[0,2\pi]
\label{variable}
\end{equation}
with initial condition $u(x,0)=1$ and the periodic boundary condition. The exact solution is given by
\begin{equation}
u(x,t) = \frac{ \sin( 2\tan^{-1} ( e^{-t}\tan(\frac{x}{2}))) }{\sin(x)}.
\end{equation}

  We report the $L^2$ and $L^\infty$ errors and the associated orders of accuracy for the proposed SLDG schemes using $\Delta t = 0.5\Delta x$ and $\Delta t = 2.5\Delta x$ in Table \ref{table:1dsin}.
 The order of convergence for the $P^1$ SLDG scheme is slightly less than
expected second order measured in both the $L^2$ and $L^\infty$ norms.
Similarly,  the order of convergence for the $P^2$ SLDG scheme is slightly
less than the expected third order.
%in the $L^2$ error and  at least second order convergence is observed in the $L^\infty$ error.
 The expected fourth order accuracy for the $P^3$ SLDG scheme is observed in both the $L^2$ and $L^\infty$ errors.
For this example, the test functions are distorted more over larger time steps. Hence, the errors from simulations with $\Delta t = 2.5\Delta x$ are observed to be slightly larger.

\begin{table}[!ht]%\small
\caption{SLDG for \eqref{variable} with $u(x,0)=1$. $T=1$.}
\vspace{0.1in}
\centering
\begin{tabular}{l ll ll ll ll}

\hline
 mesh & \multicolumn{4}{c}{$P^k$ SLDG with $\Delta t= 0.5\Delta x $}
     &\multicolumn{4}{c}{$P^k$ SLDG with $\Delta t= 2.5\Delta x $ } \\
   \cmidrule(lr){2-5} \cmidrule(lr){6-9}
{   }  &{$L^2$ error} & Order  & {$L^{\infty}$ error} & Order  &{$L^2$ error} & Order  & {$L^{\infty}$ error} & Order \\
\hline
    \multicolumn{9}{l}{ $P^1$ SLDG }
     \\

    20 &     1.28E-02 &          &    4.46E-02 &
     & 1.73E-02 &   & 6.21E-02   &  \\
    40 &     3.46E-03 &     1.89 &     1.30E-02 &     1.77
     &  1.01E-02 &     0.77 &     4.41E-02 &     0.49 \\
    80 &     9.17E-04 &     1.92 &     4.70E-03 &     1.47
     &     1.40E-03 &     2.85 &     7.26E-03 &     2.60 \\
    160 &     2.50E-04 &     1.87 &     1.40E-03 &     1.75
     &    3.24E-04 &     2.11 &     1.42E-03 &     2.35 \\
     320 &     6.47E-05 &     1.95 &     3.91E-04 &     1.84
    &    1.00E-04 &     1.69 &     3.60E-04 &     1.98  \\
    \multicolumn{9}{l}{ $P^2$ SLDG}
     \\

    20 &     1.58E-03 &          &    9.83E-03 &
     & 6.03E-03 &   &2.59E-02 &  \\
    40 &     2.45E-04 &     2.69 &     1.85E-03 &     2.41
     &8.08E-04 &     2.90 &     5.22E-03 &     2.31\\
    80 &    4.00E-05 &     2.62 &     3.94E-04 &     2.23
     &     9.09E-05 &     3.15 &     5.98E-04 &     3.13\\
    160 &     6.62E-06 &     2.59 &     9.48E-05 &     2.06
     &     1.14E-05 &     2.99 &     1.13E-04 &     2.40  \\
     320 &     1.14E-06 &     2.54 &     2.33E-05 &     2.02
    &     1.80E-06 &     2.67 &     2.56E-05 &     2.15  \\

    \multicolumn{9}{l}{ $P^3$ SLDG }
     \\

    20 &     1.23E-04 &          &    6.91E-04 &
     & 7.54E-04 &   &  2.85E-03   &  \\
    40 &     7.82E-06 &     3.97 &     4.95E-05 &     3.80
     &2.60E-05 &     4.86 &     1.61E-04 &     4.15 \\
    80 &     4.88E-07 &     4.00 &     4.04E-06 &     3.61
     &    1.52E-06 &     4.09 &     7.39E-06 &     4.45 \\
    160 &     3.42E-08 &     3.83 &     2.98E-07 &     3.76
     &     9.92E-08 &     3.94 &     7.47E-07 &     3.31 \\
     320 &     2.21E-09 &     3.95 &     2.12E-08 &     3.82
    &     7.12E-09 &     3.80 &     8.75E-08 &     3.09  \\  \hline
\end{tabular}
\label{table:1dsin}
\end{table}

\end{exa}

%%%%%%%%%%%%%%%%%%%%%%%%%%%%%%%%%%%%%%%%%%%%%%%%%%%%%%%%%%%%%%%%%%%%%%%%%%%%%%%

\begin{exa}
(2D linear equation.) Consider

\begin{equation}
\begin{cases}
u_t + u_x + u_y = 0, \ x\in[0,2\pi], \ y\in[0,2\pi] \\
u(x,y,0) = \sin(x+y)
\end{cases}
\label{linear2d}
\end{equation}
with the periodic boundary conditions in both $x$ and $y$ directions. The exact solution of \eqref{linear2d} is
\begin{equation*}
u(x,y,t) = \sin(x+y-2t).
\end{equation*}

Table \ref{table:2dlinear} summarizes the $L^2$ and $L^\infty$
errors and the associated orders of accuracy
when the  $P^k$ SLDG methods with $\Delta t = 0.5\Delta x$ and $\Delta t =2.5\Delta x$ are applied to equation \eqref{linear2d}.  Expected orders of accuracy are observed. Comparable errors are observed from simulations with these two different CFLs.

\begin{table}[!ht]%\small
\caption{SLDG for \eqref{linear2d} with $u(x,y,0)=\sin (x+y)$ at $T=\pi$.}
\vspace{0.1in}
\centering
\begin{tabular}{l ll ll ll ll}

\hline
 mesh & \multicolumn{4}{c}{$P^k$ SLDG with $\Delta t= 0.5\Delta x $}
     &\multicolumn{4}{c}{$P^k$ SLDG with $\Delta t= 2.5\Delta x $ } \\
   \cmidrule(lr){2-5} \cmidrule(lr){6-9}
{   }  &{$L^2$ error} & Order  & {$L^{\infty}$ error} & Order  &{$L^2$ error} & Order  & {$L^{\infty}$ error} & Order \\
\hline
    \multicolumn{9}{l}{ $P^1$ SLDG}
     \\

    20$\times$20 &     7.00E-03 &          &    3.24E-02 &
     & 6.86E-03 &   & 3.40E-02   &  \\
    40$\times$40 &     1.73E-03 &     2.02 &     8.40E-03 &     1.95
     &  1.72E-03 &     2.00 &     8.60E-03 &     1.98 \\
    80$\times$80 &     4.31E-04 &     2.00 &     2.13E-03 &     1.98
     &     4.30E-04 &     2.00 &     2.16E-03 &     1.99 \\
    160$\times$160 &     1.08E-04 &     2.00 &     5.38E-04 &     1.99
     &    1.08E-04 &     2.00 &     5.41E-04 &     2.00 \\

    \multicolumn{9}{l}{ $P^2$ SLDG}
     \\

    20$\times20$ &     3.50E-04 &          &    4.06E-04 &
     & 3.49E-04 &   &4.06E-04 &  \\
    40$\times$40 &     4.37E-05 &     3.00 &     5.05E-05 &     3.00
     &4.37E-05 &     3.00 &     5.05E-05 &     3.01\\
    80$\times$80 &    5.46E-06 &     3.00 &     6.31E-06 &     3.00
     &     5.46E-06 &     3.00 &     6.31E-06 &     3.00\\
    160$\times$160 &     6.83E-07 &     3.00 &     7.89E-07 &     3.00
     &     6.83E-07 &     3.00 &     7.89E-07 &     3.00   \\

  \hline
\end{tabular}
\label{table:2dlinear}
\end{table}

\end{exa}

\begin{exa}
(Rigid body rotation.) Consider
\begin{equation}
u_t - (yu)_x + (xu)_y =0, \  x\in[-2\pi,2\pi],\ y\in[-2\pi,2\pi].
\label{rigid1}
\end{equation}
with the initial condition $u(x,y,0)=\exp (-x^2-y^2)$. We  apply the proposed $P^k$ SLDG methods with $k=1,\,2$ to solve this problem. Table \ref{table:rigid} summarizes the $L^2$ and $L^\infty$ errors and the associated orders of accuracy. In the simulations, we let $\Delta t = 0.5\Delta x$ and $\Delta t = 2.5\Delta x$
and the solution is computed up to $T = 2\pi$. The second order accuracy measured with the $L^2$ error is observed as expected. However, we also observe lightly
order reduction in $L^\infty$ error.
 The third order accuracy for the $P^2$ SLDG scheme is observed  in both the $L^2$ norm and $L^\infty$ norm.

\begin{table}[!ht]%\small
\caption{SLDG for \eqref{rigid1} with $u(x,y,0)=\exp (-x^2-y^2)$ at $T=2\pi$.}
\vspace{0.1in}
\centering
\begin{tabular}{l ll ll ll ll}

\hline
 mesh & \multicolumn{4}{c}{$P^k$ SLDG with $\Delta t= 0.5\Delta x $}
     &\multicolumn{4}{c}{$P^k$ SLDG with $\Delta t= 2.5\Delta x $ } \\
   \cmidrule(lr){2-5} \cmidrule(lr){6-9}
{   }  &{$L^2$ error} & Order  & {$L^{\infty}$ error} & Order  &{$L^2$ error} & Order  & {$L^{\infty}$ error} & Order \\
\hline
    \multicolumn{9}{l}{ $P^1$ SLDG}
     \\

    20$\times$20 &     1.80E-02   &          &    1.12E-01   &
     & 1.06E-02 &   & 1.12E-01   &  \\
    40$\times$40 &     3.61E-03 &     2.32 &     3.29E-02 &     1.76
     &   3.10E-03 &     1.78 &     3.36E-02 &     1.73 \\
    80$\times$80 &     7.71E-04 &     2.23 &     1.07E-02 &     1.62
     &     6.83E-04 &     2.18 &      8.97E-03 &     1.90 \\
    160$\times$160 &     1.81E-04 &     2.09 &     3.32E-03 &     1.69
     &    1.68E-04 &     2.02 &     3.01E-03 &     1.57  \\

    \multicolumn{9}{l}{ $P^2$ SLDG}
     \\

    20$\times20$ &     1.80E-03  &          &    1.85E-02   &
     & 5.36E-03 &   &5.70E-02 &  \\
    40$\times$40 &     2.14E-04 &     3.07 &     2.84E-03 &     2.70
     &2.37E-04 &     4.50 &     2.80E-03 &     4.35\\
    80$\times$80 &    2.66E-05 &     3.01 &    4.03E-04 &     2.82
     &     2.93E-05 &     3.02 &     5.26E-04 &     2.41\\
    160$\times$160 &     3.34E-06 &     3.00 &     5.09E-05 &     2.98
     &     3.42E-06 &     3.10 &     6.28E-05 &     3.07  \\

  \hline
\end{tabular}
\label{table:rigid}
\end{table}
%%%%

\end{exa}

\begin{exa}
\label{exa1}
We numerically solve equation \eqref{rigid1} with an initial condition plotted in Figure \ref{rotation_init}, which
 consists of a slotted disk, a cone as well as a
smooth hump, similar to the one used in \cite{leveque1996high} for comparison purposes.
The numerical solutions computed
by the $P^k$ SLDG methods with $k=1,\,2$ after one full rotation are plotted in Figure \ref{rotation}. Note that the schemes are able to effectively resolve the complex solution structures, while some mild oscillations appear in the vicinity of discontinuities. Once
the WENO limiter and the BP filter are applied, the spurious oscillations are suppressed and the positivity of the numerical solution is guaranteed. In addition, we observe that the $P^2$ SLDG scheme is able to better resolve solution structures compared with the $P^1$ counterpart.

\begin{figure}[h!]
\centering                              %ʹ²åͼ¾ÓÖÐ
\includegraphics[height=65mm]{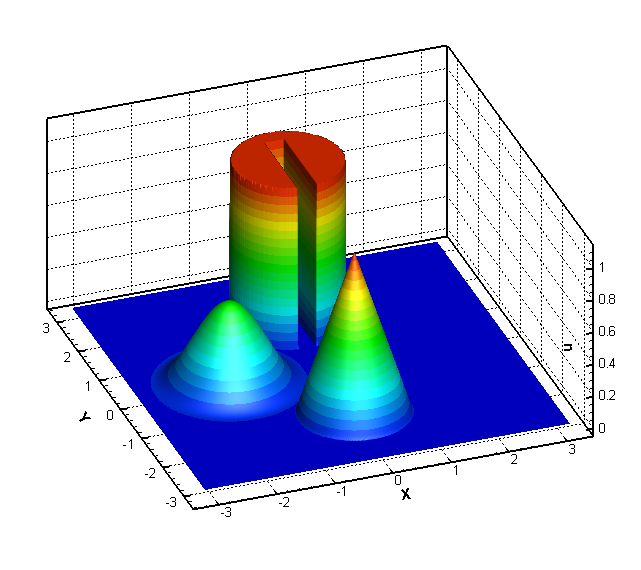}
\caption{Plots of the initial profile.  The mesh is 400$\times$400.}
\label{rotation_init}
\end{figure}

\begin{figure}[h!]
\centering                              %ʹ²åͼ¾ÓÖÐ
\subfigure[$P^1$ SLDG]{
\includegraphics[height=65mm]{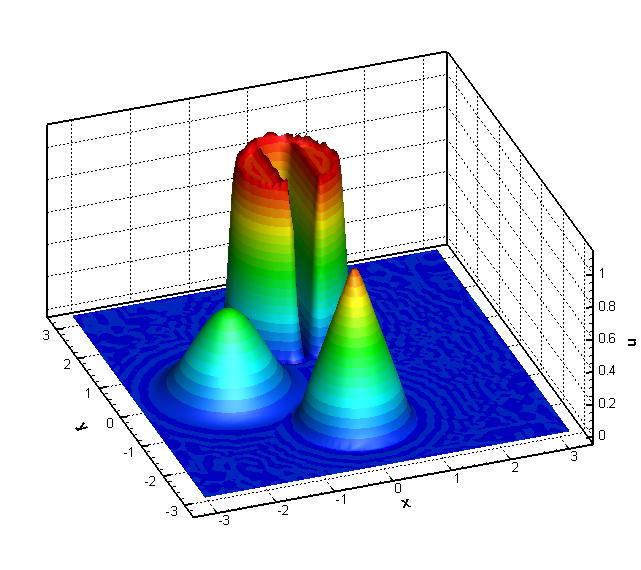} }
\subfigure[$P^1$ SLDG+WENO+BP]{
\includegraphics[height=65mm]{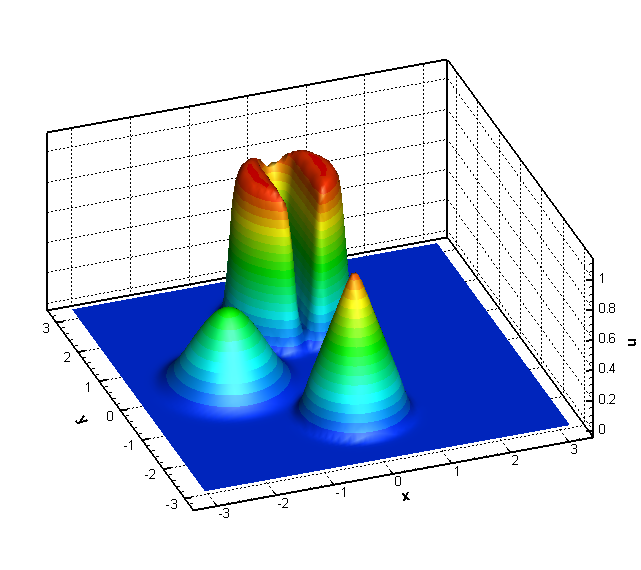} }
\subfigure[$P^2$ SLDG]{
\includegraphics[height=65mm]{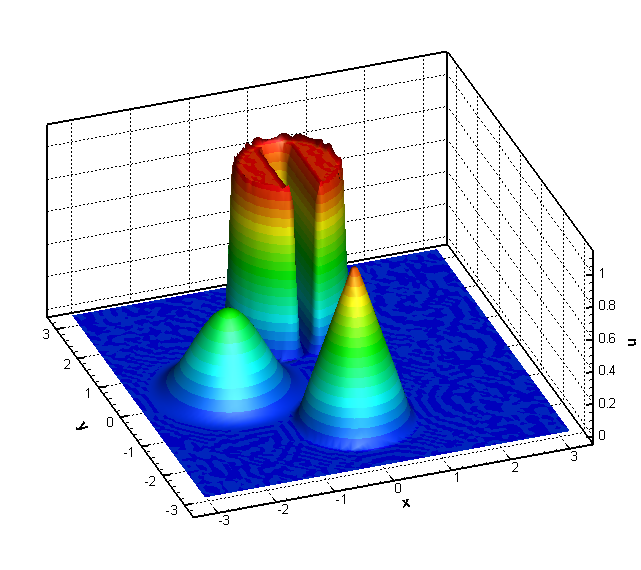} }
\subfigure[$P^2$ SLDG+WENO+BP]{
\includegraphics[height=65mm]{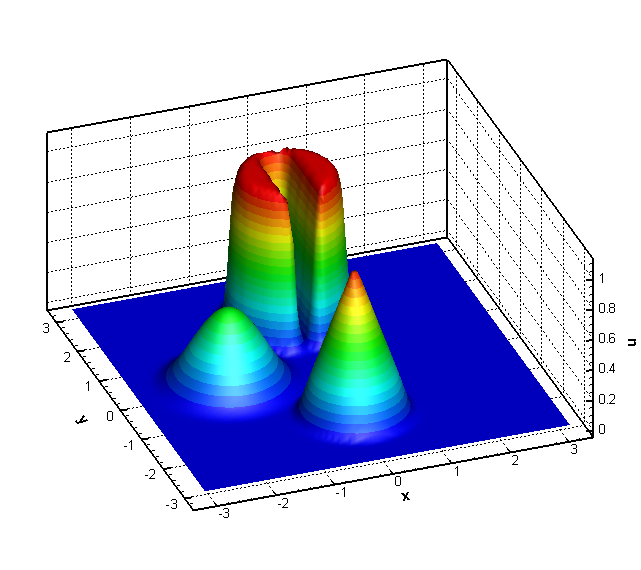} }

\caption{Plots of the numerical solutions of SLDG methods for equation \eqref{rigid1} with initial data Figure \ref{rotation_init};  $T = 2\pi$; The numerical mesh has a resolution of 80$\times$80. $\Delta t = 2.5\Delta x$. Top left: $P^1$ SLDG; Top right: $P^1$ SLDG with WENO limiter and BP filter;
Bottom left: $P^2$ SLDG; Bottom right: $P^2$ SLDG with WENO limiter and BP filter.}
\label{rotation}
\end{figure}

\end{exa}

\begin{exa}(Swirling deformation flow.)
We consider solving
\begin{equation}
u_t - \left( \cos^2 \left(\frac{x}{2} \right)\sin(y) g(t) u \right)_x
+ \left(  \sin(x) \cos^2 \left(\frac{y}{2}\right) g(t) u \right)_y =0,
\
x\in[-\pi,\pi],\ y\in [-\pi,\pi],
\label{swirling}
\end{equation}
 where $g(t) = \cos\left( \frac{\pi t}{T} \right)\pi$. The initial condition is set to be the following smooth cosine bells (with $C^5$ smoothness),
\begin{equation}
u(x,y,0) =
\begin{cases}
r_0^b \cos^6 \left(  \frac{r^b(\mathbf{x}) \pi }{2r_0^b}  \right), & \text{if}\  r^b(\mathbf{x}) <r_0^b,\\
0,  & \text{ otherwise},
\end{cases}
\label{cosine_bell}
\end{equation}
where $r_0^b = 0.3\pi$, and $r^b(\mathbf{x})=\sqrt{  (x-x_0^b)^2 + (y-y_0^b)^2 }$ denotes the distance between $(x,y)$ and the center of the cosine bell $(x_0^b, y_0^b) = ( 0.3\pi ,0 )$.
In  Table \ref{table:df_order}, we summarize the
convergence study for $P^1$ and $P^2$ SLDG as well as $P^2$ SLDG-QC methods in terms of the $L^2$ and $L^\infty$ errors. We observe the second order convergence for the $P^1$ and $P^2$ SLDG schemes  measured by both errors. Furthermore, the third order convergence for the $P^2$ SLDG-QC scheme  measured by the $L^2$ error  is observed. Half order reduction is observed for the $L^\infty$ error. In addition, the magnitude of the errors for the $P^2$ SLDG-QC method is much smaller than both $P^1$ and $P^2$ SLDG methods. The observation justifies the use of the  quadratic-curved quadrilateral approximation in the SLDG formulation.

\begin{table}[!ht]%\small
\caption{SLDG for \eqref{swirling} with \eqref{cosine_bell} at $T=1.5$.}
\vspace{0.1in}
\centering
\begin{tabular}{l ll ll ll ll}

\hline
 mesh & \multicolumn{4}{c}{$P^k$ SLDG with $\Delta t= 0.5\Delta x $}
     &\multicolumn{4}{c}{$P^k$ SLDG with $\Delta t= 2.5\Delta x $ } \\
   \cmidrule(lr){2-5} \cmidrule(lr){6-9}
{   }  &{$L^2$ error} & Order  & {$L^{\infty}$ error} & Order  &{$L^2$ error} & Order  & {$L^{\infty}$ error} & Order \\
\hline
    \multicolumn{9}{l}{ $P^1$ SLDG}
     \\

    20$\times$20 &     1.25E-02 &       &    2.03E-01 &
   &    8.59E-03 &     &     2.45E-01 &     \\
    40$\times$40 &    2.92E-03 &     2.10 &     6.92E-02 &     1.55
  & 2.14E-03 &     2.00 &     6.74E-02 &     1.87  \\
    80$\times$80 &    5.96E-04 &     2.29 &     1.48E-02 &     2.22
  &     5.42E-04 &     1.98 &     2.19E-02 &     1.62  \\
   160$\times$160 &     1.30E-04 &     2.20 &     4.32E-03 &     1.78
   &    1.33E-04 &     2.02 &     6.16E-03 &     1.83 \\

    \multicolumn{9}{l}{ $P^2$ SLDG}
     \\

    20$\times$20 &    3.22E-03 &        &     8.97E-02 &
   &     9.37E-03 &       &      2.70E-01 &       \\
    40$\times$40 &    6.58E-04 &     2.29 &     3.44E-02 &     1.38
  & 2.87E-03 &     1.71 &     1.28E-01 &     1.08 \\
    80$\times$80 &    1.42E-04 &     2.22 &     8.37E-03 &     2.04
  &   6.92E-04 &     2.05 &     4.19E-02 &     1.61   \\
   160$\times$160 &     3.15E-05 &     2.17 &     2.47E-03 &     1.76
   &    1.89E-04 &     1.87 &     1.20E-02 &     1.80  \\

      \multicolumn{9}{l}{ $P^2$ SLDG-QC}
     \\

    20$\times$20 &     2.61E-03 & &     6.27E-02 &
   &     5.29E-03 & &     1.51E-01 &       \\
    40$\times$40 &     3.15E-04 &     3.05 &     1.20E-02 &     2.38
  & 7.78E-04 &     2.77 &     2.63E-02 &     2.52  \\
    80$\times$80 &   3.81E-05 &     3.05 &     1.78E-03 &     2.76
  &     1.04E-04 &     2.90 &     4.59E-03 &     2.52   \\
   160$\times$160 &     4.91E-06 &     2.96 &     2.74E-04 &     2.70
   &   1.47E-05 &     2.83 &     7.85E-04 &     2.55 \\ \hline
\end{tabular}
\label{table:df_order}
\end{table}

\end{exa}

\begin{exa}
 We numerically solve equation \eqref{swirling} with $g(t) = \cos\left( \frac{\pi t}{T} \right)\pi$ and the same initial condition as in Example \ref{exa1}.
Numerical solutions for the SLDG methods are
plotted in Figure \ref{df1}.
To better compare performance of the schemes with different configurations, we plot
 1D cuts of the numerical solutions along with the exact solution in Figure \ref{df2}.
 In Figure \ref{df1} and Figure \ref{df2}, one can observe that,  without using limiters,
 some mild spurious oscillations appear, which can be removed by further coupling the WENO limiter. Again, positivity of the numerical solution is guaranteed when the BP filter is applied.
 Note that, for this test, the performance between the $P^2$ SLDG-QC and $P^2$ SLDG schemes is qualitatively comparable, which is slightly better compared with the  $P^1$ SLDG scheme.
   Figure \ref{df3} shows the contour plots of the numerical
solution for SLDG methods at the final integration time 0.75,
when the solution is quite deformed. Figure \ref{df4} presents the time evolution of the relative error in the total mass, which demonstrates the mass conservation property of the proposed SLDG schemes.

\begin{figure}[h!]
\centering                              %ʹ²åͼ¾ÓÖÐ
\subfigure[$P^1$ SLDG]{
\includegraphics[height=65mm]{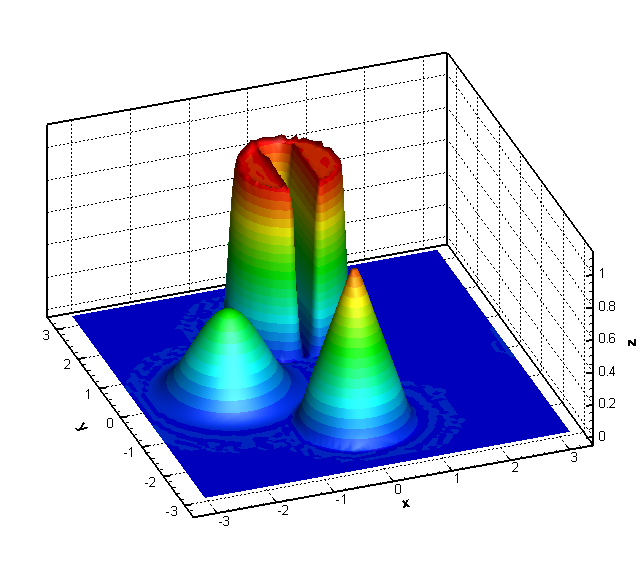} }
\subfigure[$P^1$ SLDG+WENO+BP]{
\includegraphics[height=65mm]{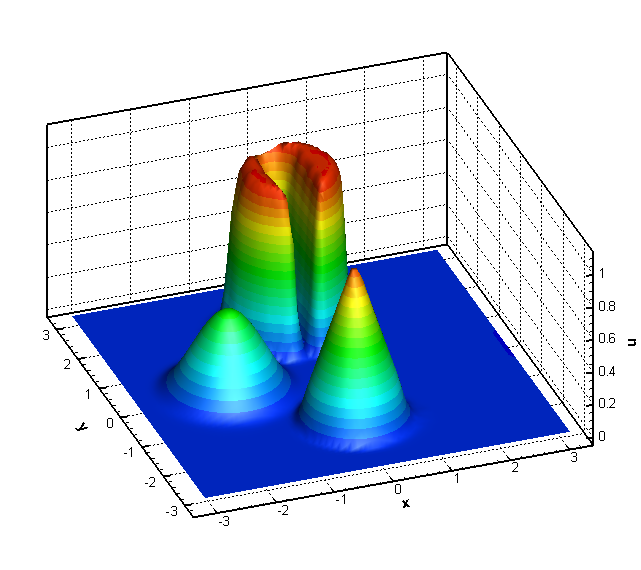} }
\subfigure[$P^2$ SLDG]{
\includegraphics[height=65mm]{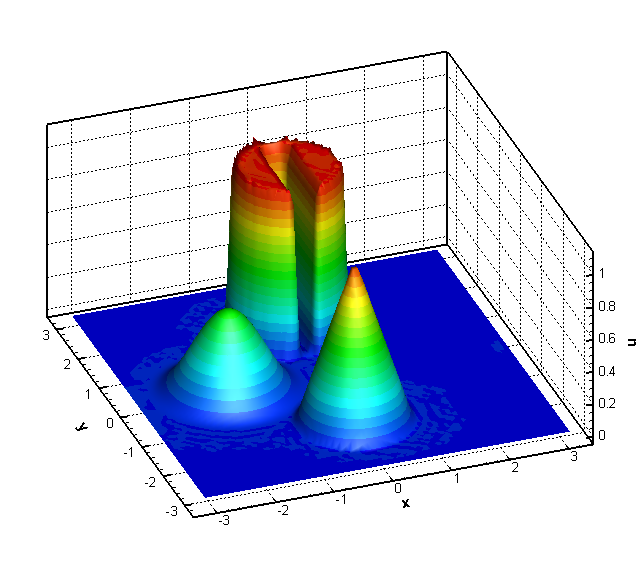} }
\subfigure[$P^2$ SLDG+WENO+BP]{
\includegraphics[height=65mm]{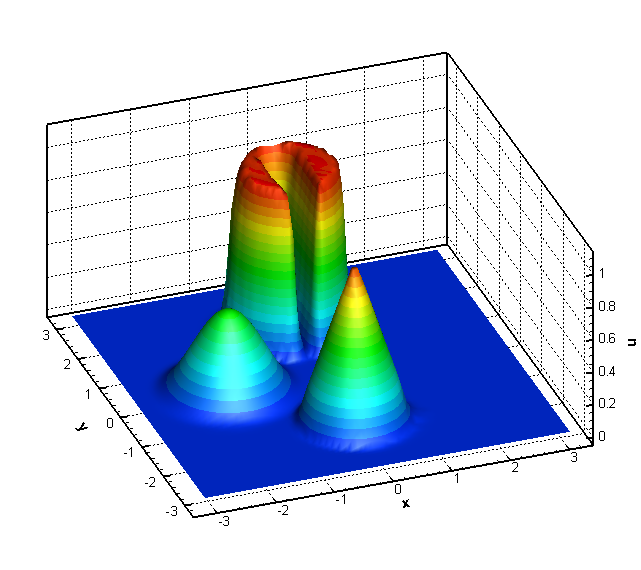} }
\subfigure[$P^2$ SLDG-QC]{
\includegraphics[height=65mm]{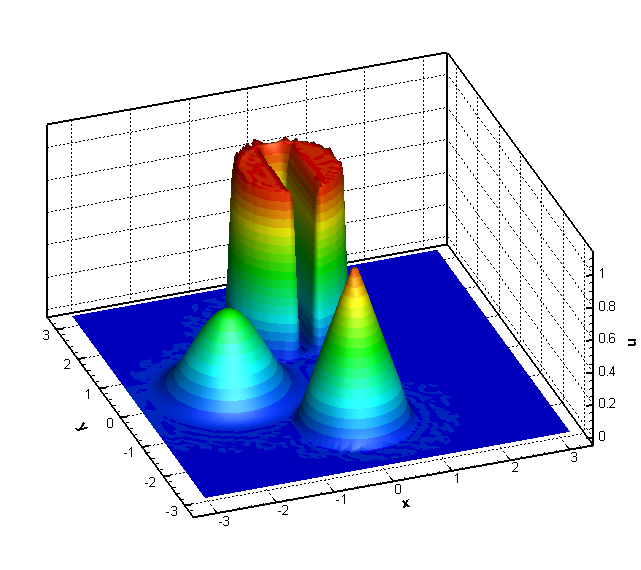}}
\subfigure[$P^2$ SLDG-QC+WENO+BP]{
\includegraphics[height=65mm]{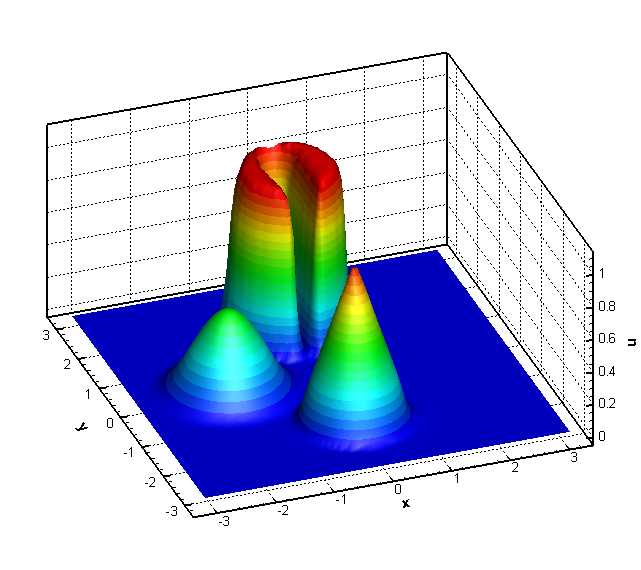}  }
\caption{
Plots of the numerical solution of SLDG for equation \eqref{swirling} with $g(t) = \cos\left( \frac{\pi t}{T} \right)\pi$.;  $T = 1.5$ and the final integration time 1.5; The numerical mesh has a resolution of 80$\times$80. $\Delta t = 2.5\Delta x$. }
\label{df1}
\end{figure}

\begin{figure}[h!]
\centering
\subfigure[SLDG]{
\includegraphics[height=75mm]{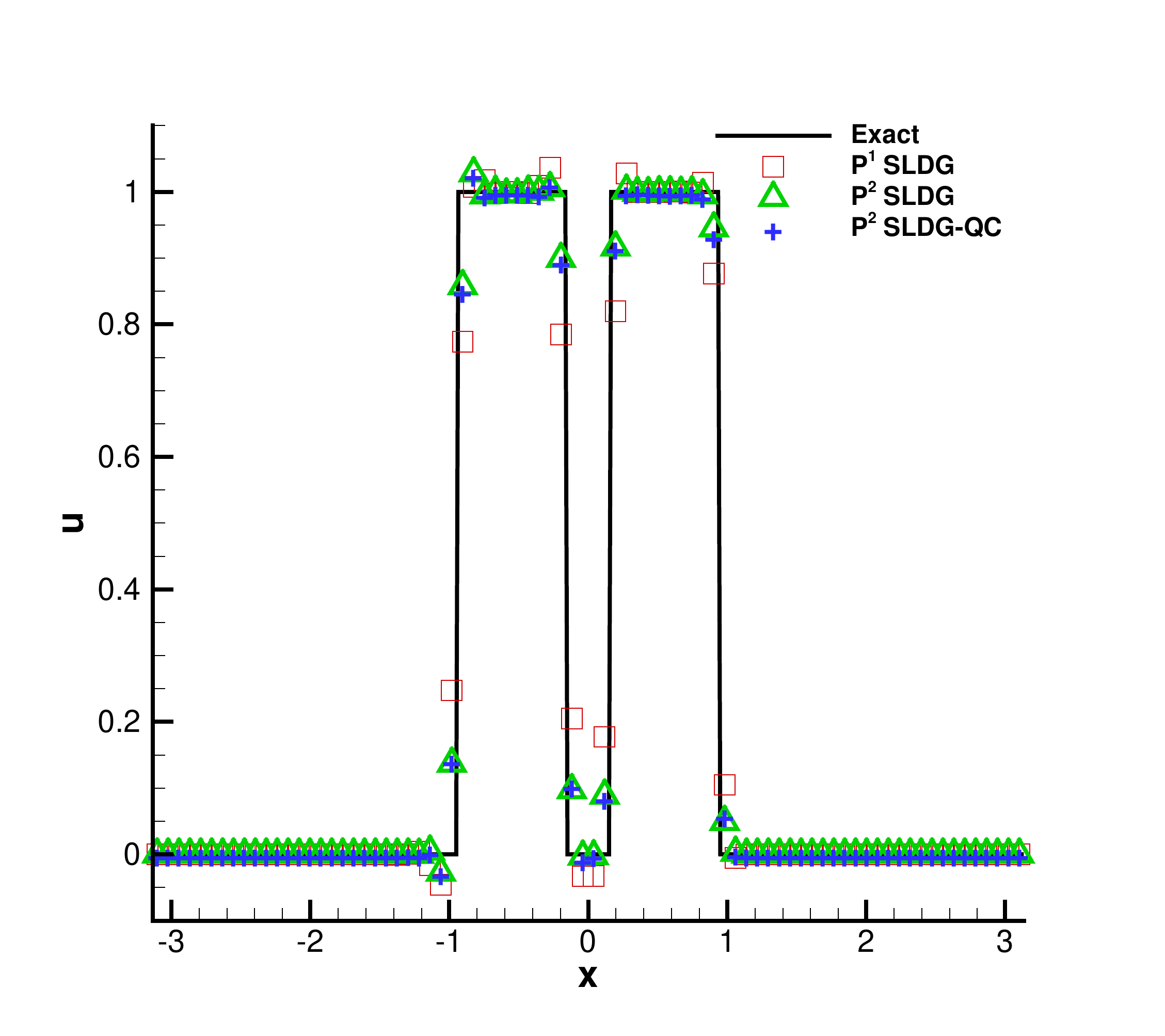} }
\subfigure[SLDG+WENO+BP]{
\includegraphics[height=75mm]{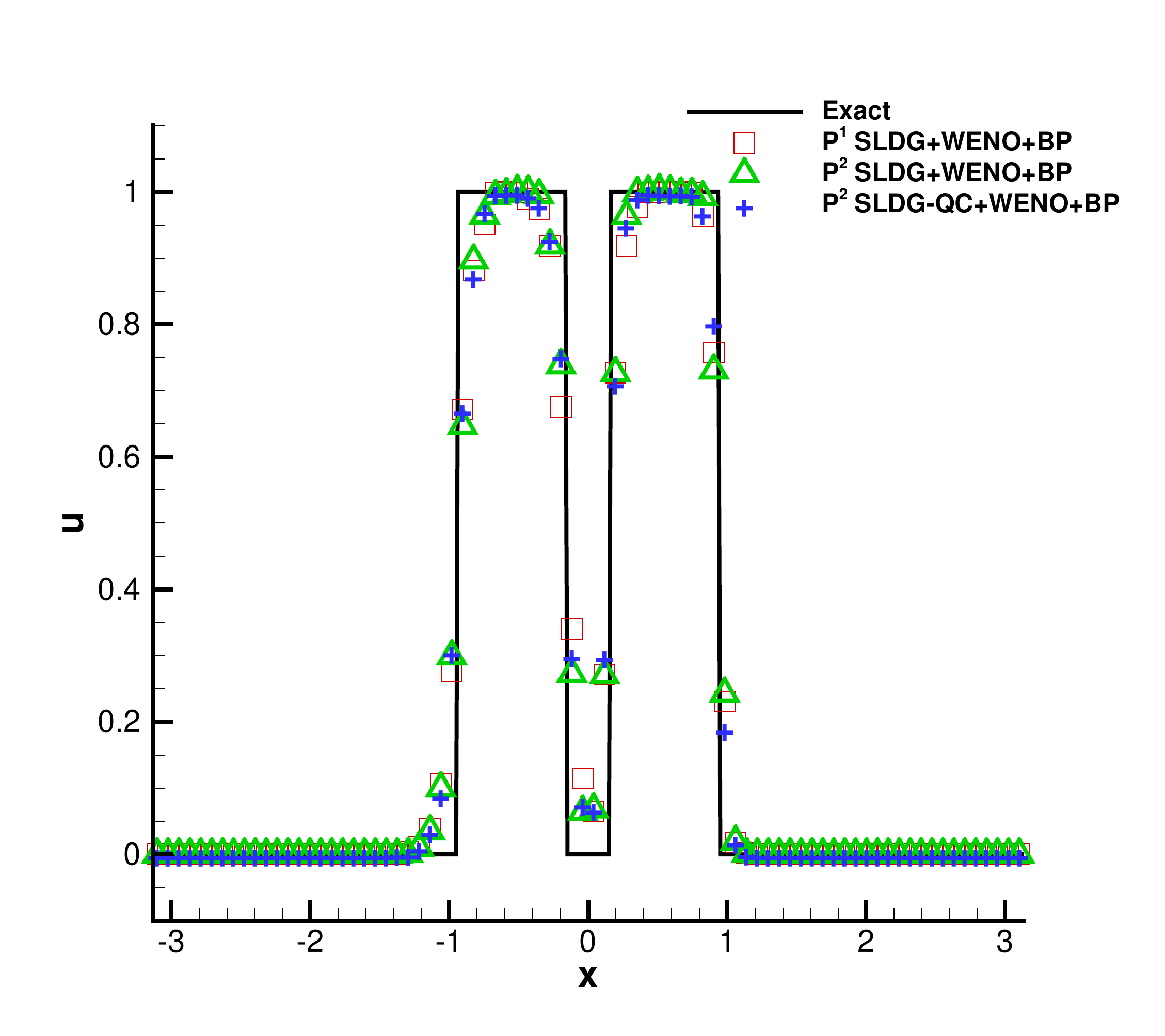} }
\caption{Plots of the 1D cuts of the numerical solution for equation \eqref{swirling} with $g(t) = \cos\left( \frac{\pi t}{T} \right)\pi$. at $y=\frac{\pi}{2}$. The solid line depicts the exact solution. The numerical mesh has a resolution of $80\times80$.}
\label{df2}
\end{figure}

\begin{figure}[h!]
\centering                              %ʹ²åͼ¾ÓÖÐ
\subfigure[Contour of $P^1$ SLDG]{
\includegraphics[height=65mm]{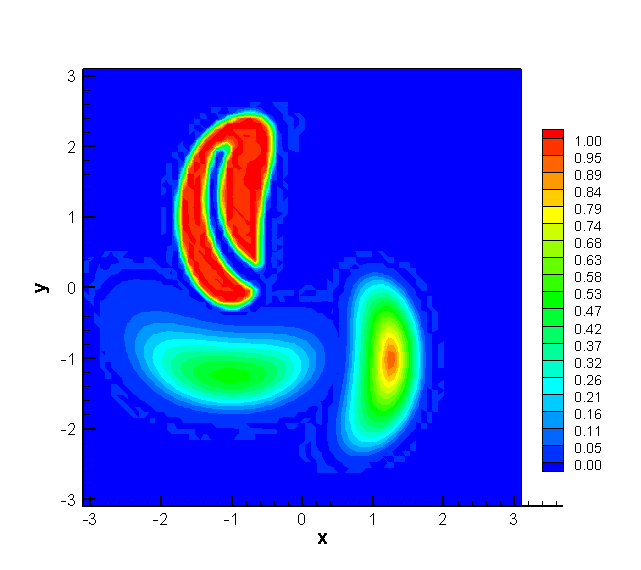} }
\subfigure[Contour of $P^1$ SLDG+WENO+BP ]{
\includegraphics[height=65mm]{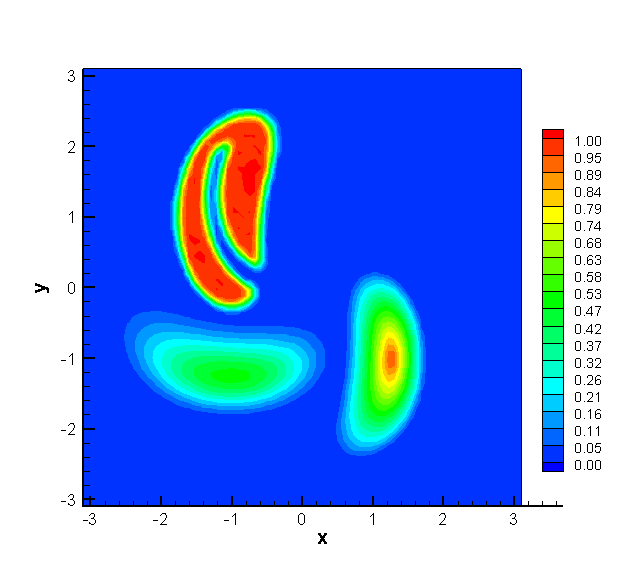} }
\subfigure[Contour of $P^2$ SLDG]{
\includegraphics[height=65mm]{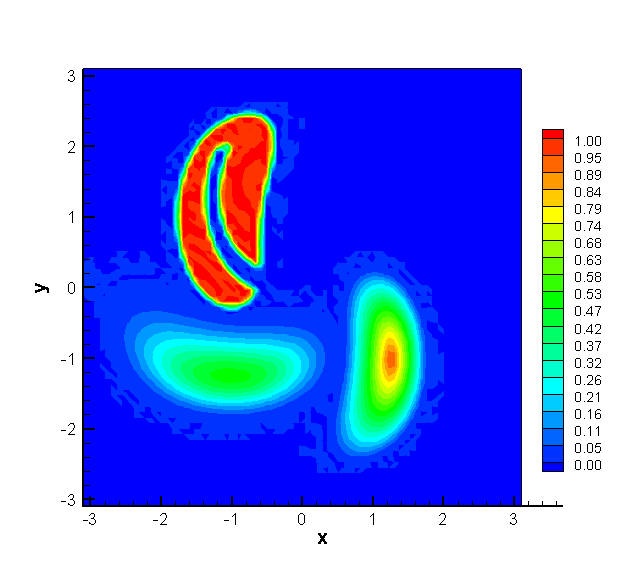} }
\subfigure[Contour of $P^2$ SLDG+WENO+BP]{
\includegraphics[height=65mm]{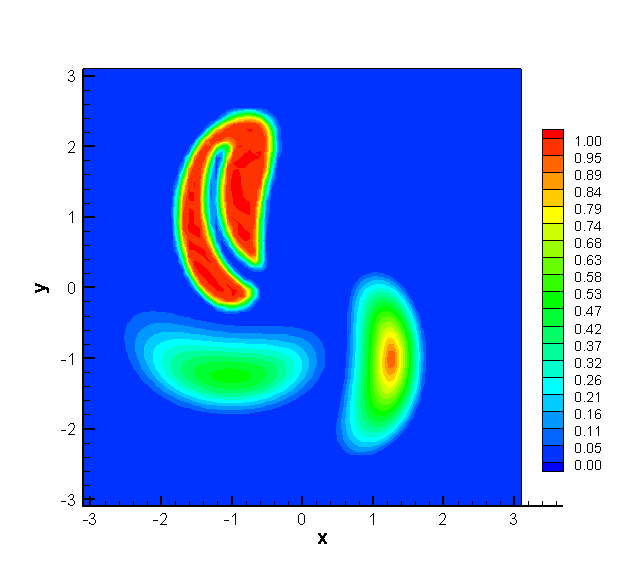} }
\subfigure[Contour of $P^2$ SLDG-QC]{
\includegraphics[height=65mm]{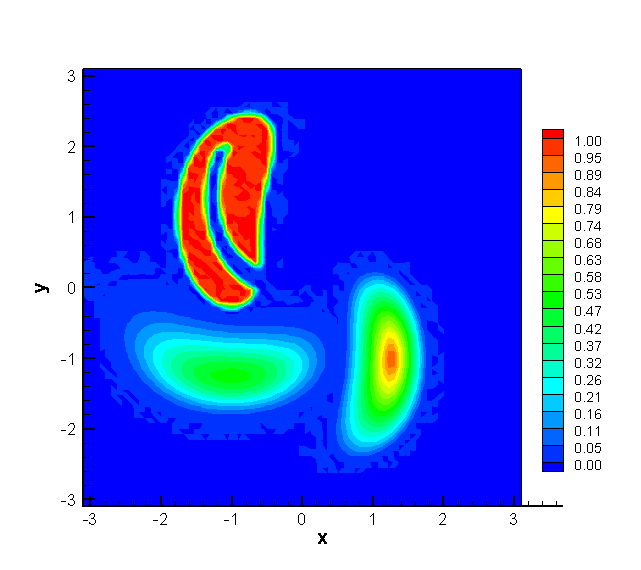}}
\subfigure[Contour of $P^2$ SLDG-QC+WENO+BP]{
\includegraphics[height=65mm]{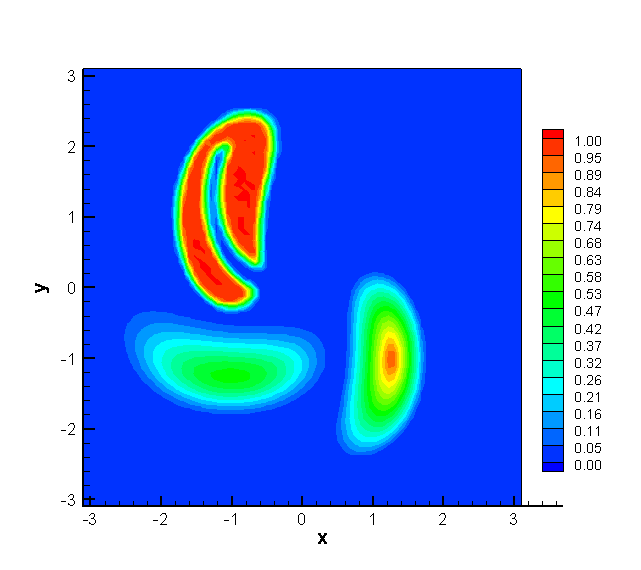}}
\caption{
Plots of the numerical solution of SLDG for equation \eqref{swirling} with $g(t) = \cos\left( \frac{\pi t}{T} \right)\pi$.;  $T = 1.5$ and the final integration time 0.75; The numerical mesh has a resolution of 80$\times$80. $\Delta t = 2.5\Delta x$.}
\label{df3}
\end{figure}

\begin{figure}[h!]
\centering                              %ʹ²åͼ¾ÓÖÐ
\subfigure[$P^1$ SLDG]{
\includegraphics[height=65mm]{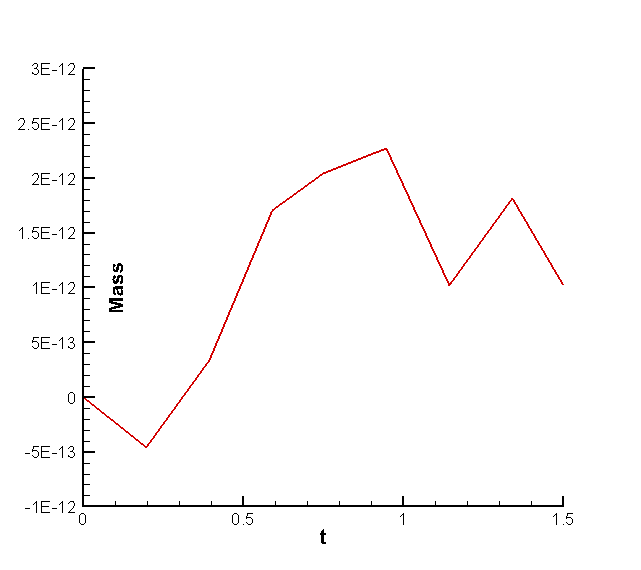} }
\subfigure[$P^2$ SLDG]{
\includegraphics[height=65mm]{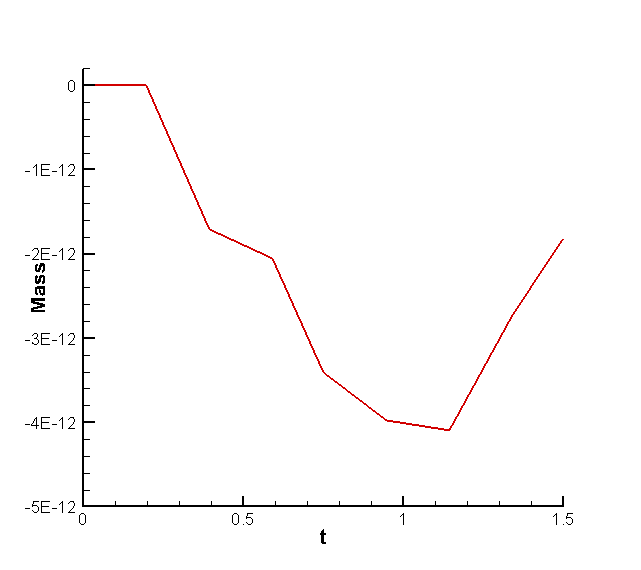} }
\subfigure[$P^2$ SLDG-QC]{
\includegraphics[height=65mm]{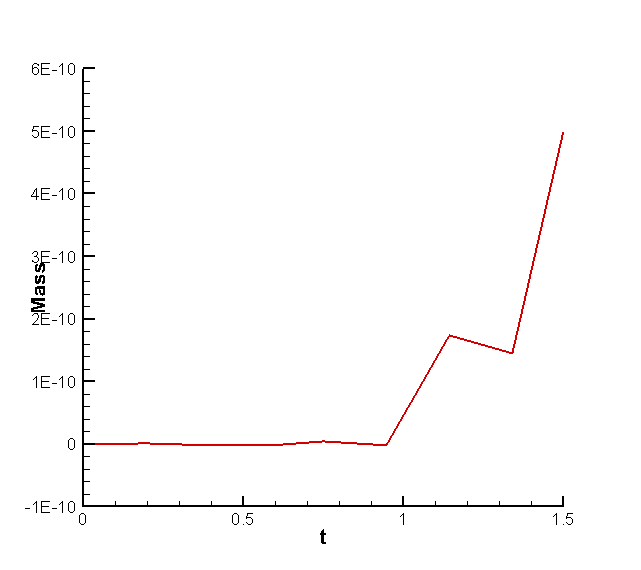} }

\caption{
The relative mass error for SLDG.
The numerical mesh has a resolution of 80$\times$80. $\Delta t = 2.5\Delta x$.
}
\label{df4}
\end{figure}

\end{exa}

\begin{exa}
 In this example, we still consider the  swirling deformation flow \eqref{swirling} but with  $g(t) = 1$.
 The initial condition is set to be
\begin{equation}
u(x,y,0) =
\begin{cases}
1, & \text{if} \ \sqrt{ (x-\pi)^2 + (y-\pi)^2 } \leq \frac{8\pi}{5},\\
0, & \text{otherwise}.
\end{cases}
\end{equation}

 The numerical solutions at $T=5\pi$ of SLDG methods are plotted in Figure \ref{swirling_g1}.
 The results with excellent resolution are observed. Again, the WENO limiter and BP filter effectively remove the undesired oscillations and keep the numerical solution nonnegative.

\begin{figure}[h!]
\centering                              %ʹ²åͼ¾ÓÖÐ
\subfigure[$P^1$ SLDG]{
\includegraphics[height=65mm]{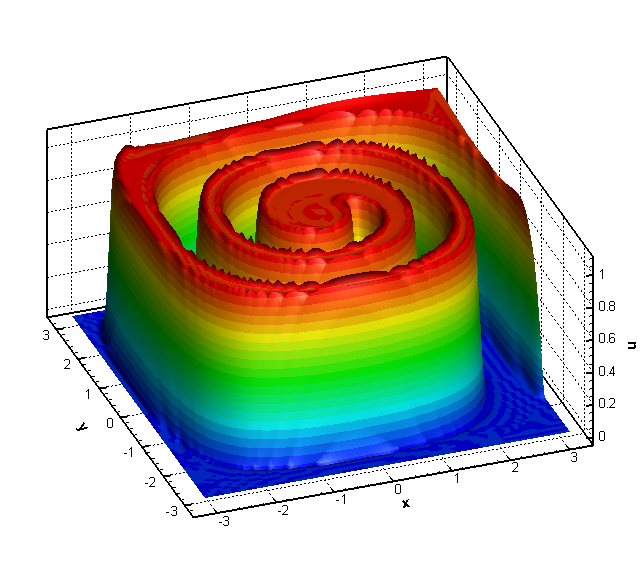}  }
\subfigure[$P^1$ SLDG+WENO+BP]{
\includegraphics[height=65mm]{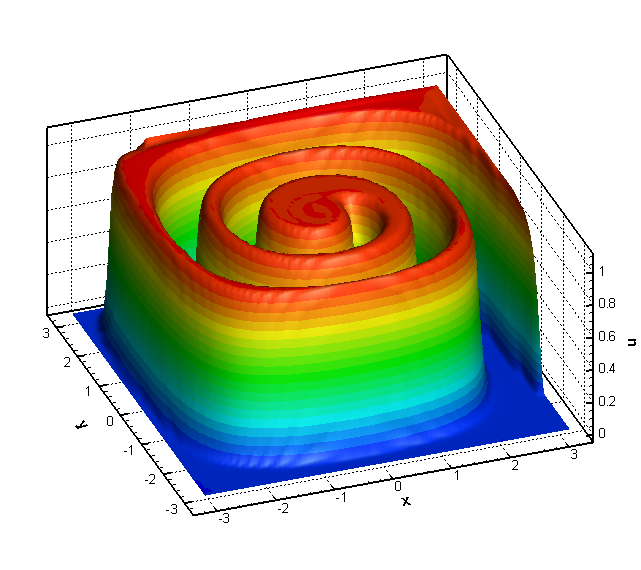}  }
\subfigure[$P^2$ SLDG]{
\includegraphics[height=65mm]{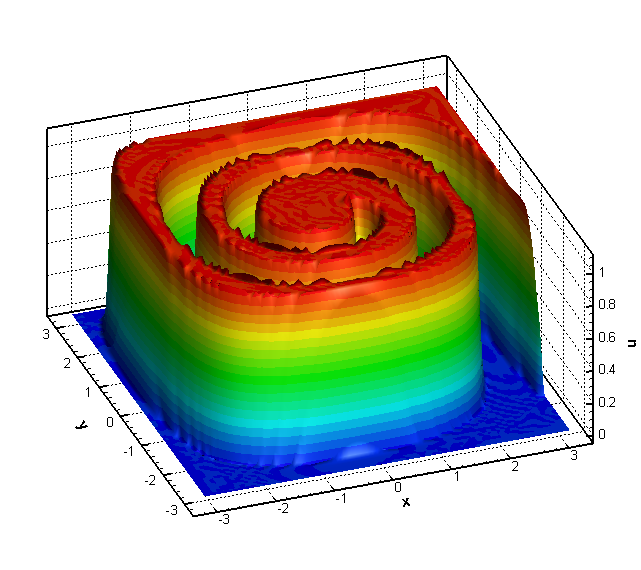} }
\subfigure[$P^2$ SLDG+WENO+BP]{
\includegraphics[height=65mm]{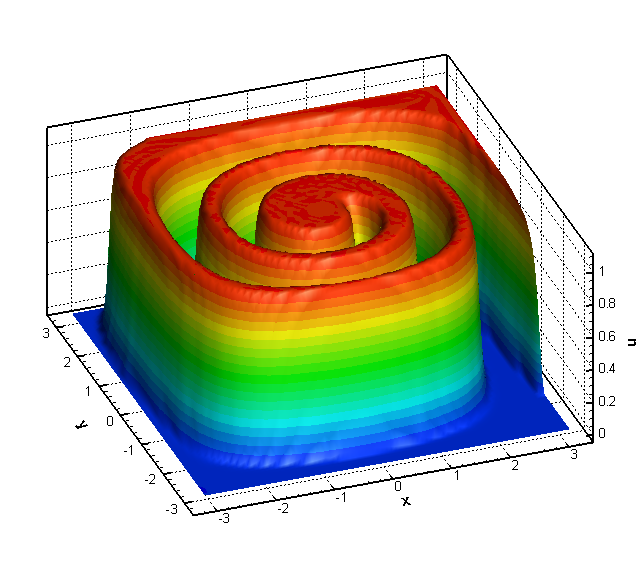}  }
\subfigure[$P^2$ SLDG-QC]{
\includegraphics[height=65mm]{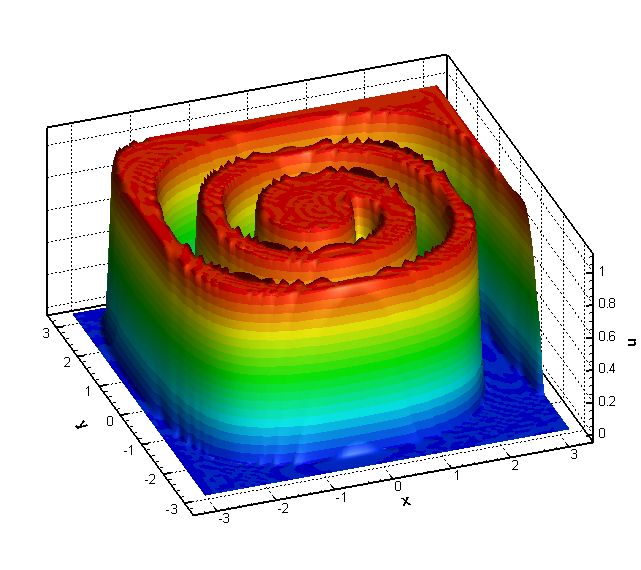} }
\subfigure[$P^2$ SLDG-QC+WENO+BP]{
\includegraphics[height=65mm]{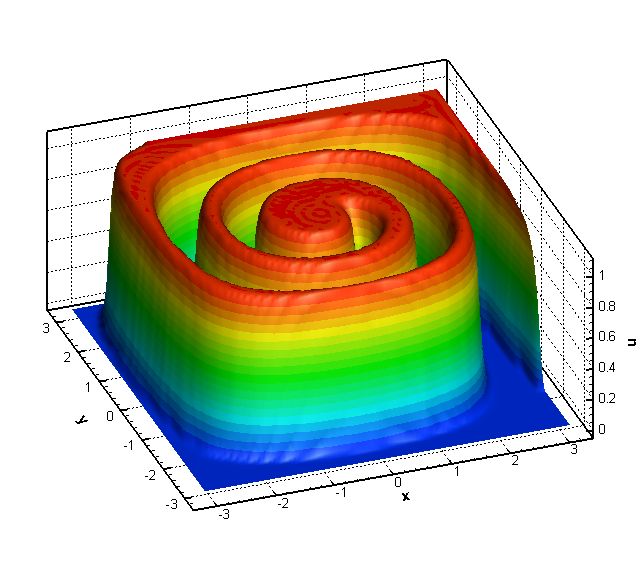} }
\caption{ Plots of the numerical solution with SLDG for equation \eqref{swirling} with $g(t)=1$. Mesh size: $80\times80$. Final integration
time $5\pi$. $\Delta t = 2.5\Delta x$. }
\label{swirling_g1}
\end{figure}

\label{example:df}
\end{exa}

\begin{exa}
We consider a 2D deformation flow in \cite{blossey2008selective}. The initial circular distribution deforms into crescent shape as it moves in the domain, and returns to the initial position when the flow reverses. The velocity field is defined on unit square $[0,1]^2 $ as
\begin{equation}
a(x,y,t)  = a_\theta(r,t) \sin (\theta),\
b(x,y,t)  = a_\theta(r,t) \cos (\theta),
\end{equation}
where $r=\sqrt{ (x-0.5)^2  + (y-0.5)^2 }$, $\theta=\tan^{-1}[ (y-0.5)/(x-0.5) ]$, and
\begin{equation*}
a_{\theta}(r,t) = \frac{ 4\pi r }{T} \left[  1 - \cos\left( \frac{2\pi r}{T}  \right) \frac{1-(4r)^6}{1+(4r)^6}  \right].
\end{equation*}
The initial condition is given by
\begin{equation*}
u(x,y,t=0) =
\begin{cases}
 \left(  \frac{1+\cos(\pi\widetilde{r} ) }{2} \right)^2 & \text{if}\ \widetilde{r}\leq 1\\
 0 & \text{if} \ \widetilde{r} >1,
\end{cases}
\end{equation*}
where $\widetilde{r} = 5\sqrt{ (x-0.3)^2 + (y-0.5)^2 }$.
Though the velocity field is complicated, the analytical solution is known at the final time $t=T$ and equals to the initial state, therefore, the error measures can be computed at time $T$.

For the simulation, the SLDG schemes with mesh 60$\times$60 and $\Delta t = 1.5\Delta x$ are used.
In Figure \ref{df_p1},  Figure \ref{df_p2QC} and Figure \ref{df_p2QC_pp}, we report the contour plots of the numerical solutions of the $P^1$ SLDG and $P^2$ SLDG-QC schemes at (a) $t=T$, (b) $t=T/4$, (c) $t=T/2$ and (d) $t=3T/4$, without and with the BP filter, respectively. It is observed that the SLDG schemes are able to capture fine features of the solution (Figure \ref{df_p1} and Figure \ref{df_p2QC}) even without using the BP filter. Compared to the solution by the $P^1$ SLDG scheme reported in Figure \ref{df_p1},  the $P^2$ SLDG-QC scheme generates the numerical solution with better resolution, see Figure \ref{df_p2QC}. Note that the BP filter ensures the positivity of numerical solution and also helps to get rid of unphysical oscillations. In Table \ref{table:rd_order}, we observe that the magnitude of error by the $P^2$ SLDG scheme is smaller than that by the $P^1$ SLDG scheme but larger than that by the $P^2$ SLDG-QC scheme, as expected. On the other hand,  the $P^2$ SLDG-QC scheme requires the most CPU time among the three schemes in which case the search algorithm is the most complicated.

\begin{figure}[h!]
\centering                              %ʹ²åͼ¾ÓÖÐ
\subfigure[Initial (Final) Field ($P^1$ SLDG)]{
\includegraphics[height=65mm]{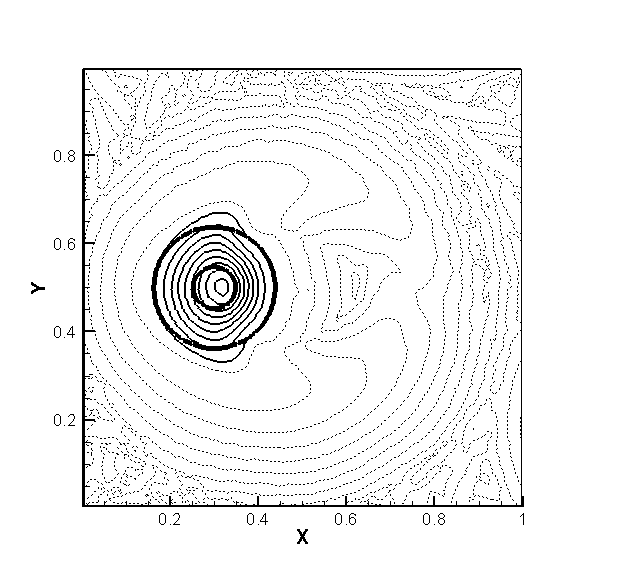} }
\subfigure[Field at T/4 ($P^1$ SLDG)]{
\includegraphics[height=65mm]{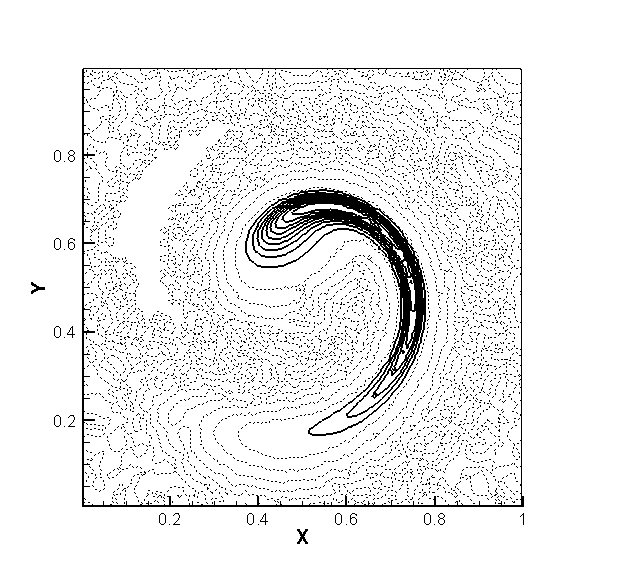} }
\subfigure[Field at T/2 ($P^1$ SLDG)]{
\includegraphics[height=65mm]{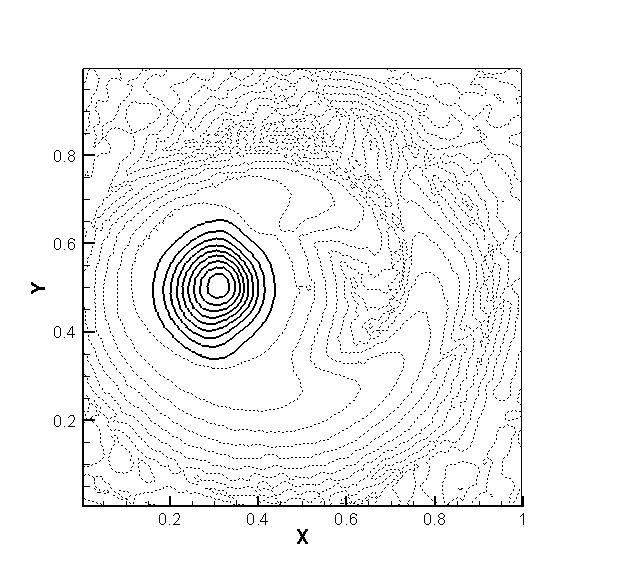} }
\subfigure[Field at 3T/4 ($P^1$ SLDG)]{
\includegraphics[height=65mm]{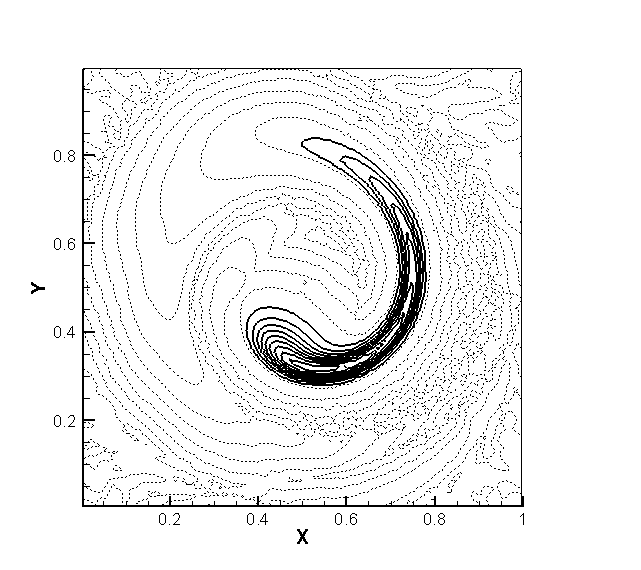} }
\caption{
2-D deformation flow. The $P^1$ SLDG scheme without the BP filter is used.
The numerical mesh has a resolution of 60$\times$60. $\Delta t = 1.5\Delta x$.
Panels (b), (c), (d) and (a) shows the deformation of the initial distribution during the simulation at time $T/4$, $T/2$, $3T/4$ and $T$, respectively. The contours are plotted in the range from -0.05 to 0.95 with increment of 0.1, and an additional contour at 0 is added (see dashed line). Numerical oscillations appear and negative numerical solution is observed. Thick contours are the highlighted exact (initial) solution for the contour values 0.05 and 0.75.
}
\label{df_p1}
\end{figure}

\begin{figure}[h!]
\centering                              %ʹ²åͼ¾ÓÖÐ
\subfigure[Initial (Final) Field ($P^2$ SLDG-QC)]{
\includegraphics[height=65mm]{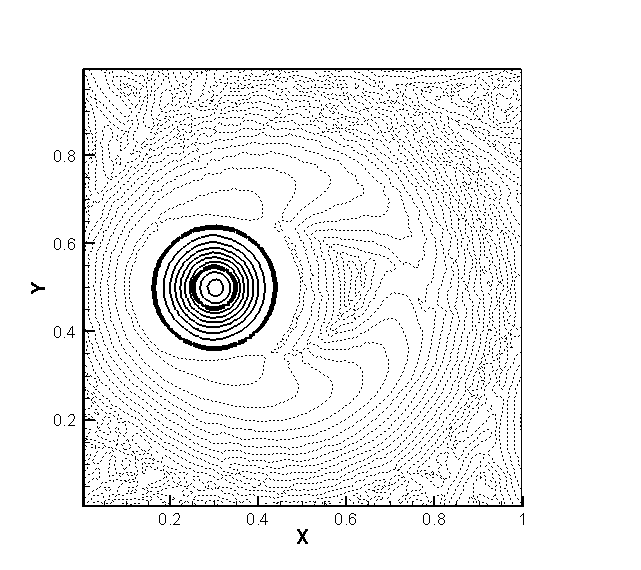} }
\subfigure[Field at T/4 ($P^2$ SLDG-QC)]{
\includegraphics[height=65mm]{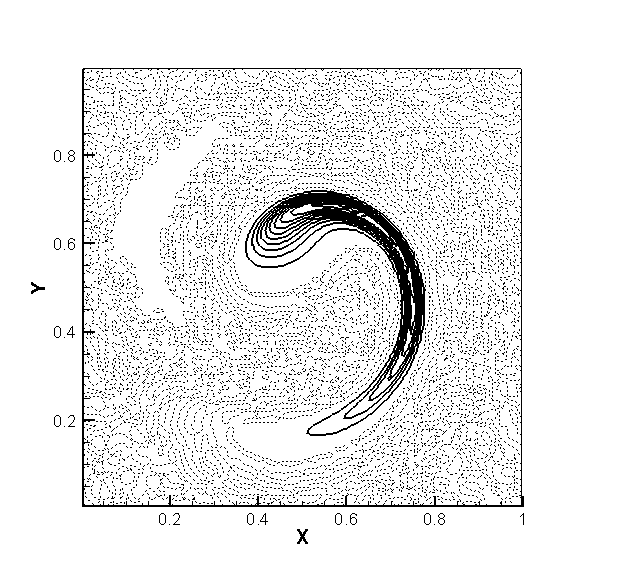} }
\subfigure[Field at T/2 ($P^2$ SLDG-QC)]{
\includegraphics[height=65mm]{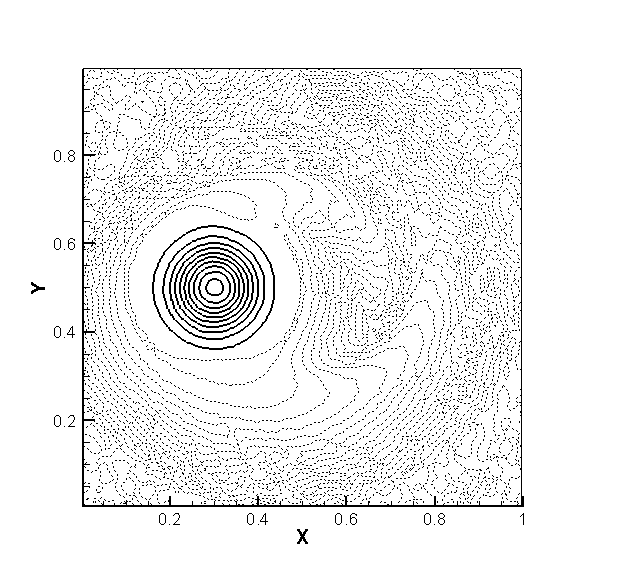} }
\subfigure[Field at 3T/4 ($P^2$ SLDG-QC)]{
\includegraphics[height=65mm]{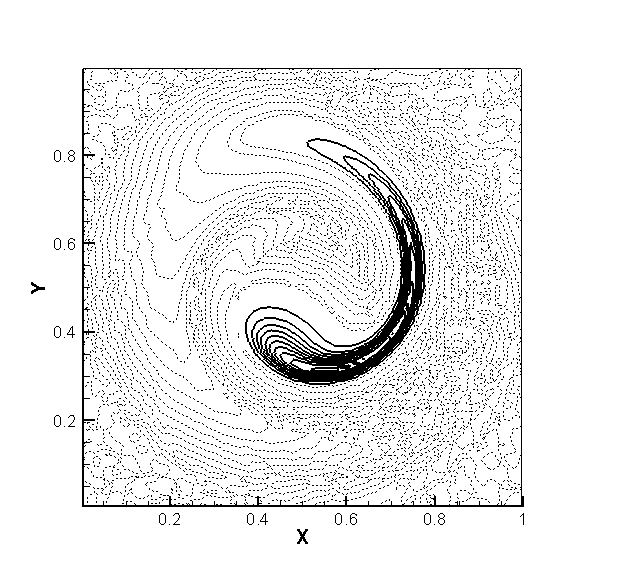} }
\caption{
2-D deformation flow. The $P^2$ SLDG scheme with quadratic-curved upstream cell without the BP filter is used.
The numerical mesh has a resolution of 60$\times$60. $\Delta t = 1.5\Delta x$.
Panels (b), (c), (d) and (a) shows the deformation of the initial distribution during the simulation at time $T/4$, $T/2$, $3T/4$ and $T$, respectively. The contours are plotted in the range from -0.05 to 0.95 with increment of 0.1, and an additional contour at 0 is added (see dashed line). Numerical oscillations appear and negative numerical solution is observed. Thick contours are the highlighted exact (initial) solution for the contour values 0.05 and 0.75.
}
\label{df_p2QC}
\end{figure}

\begin{figure}[h!]
\centering                              %ʹ²åͼ¾ÓÖÐ
\subfigure[Initial (Final) Field ($P^2$ SLDG-QC+BP)]{
\includegraphics[height=65mm]{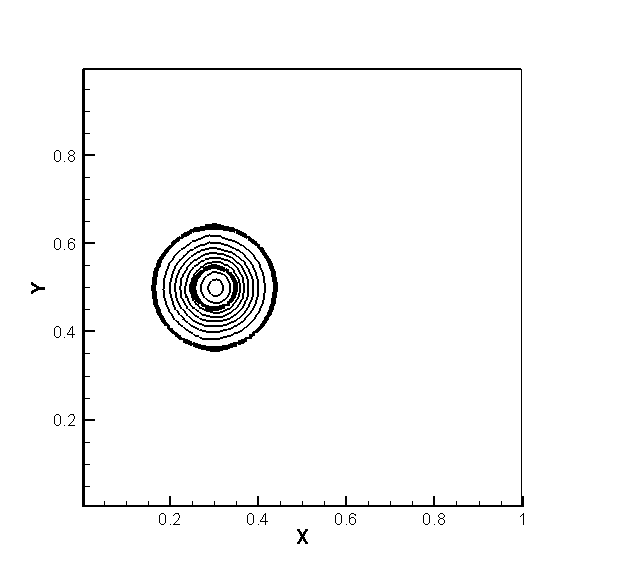} }
\subfigure[Field at T/4 ($P^2$ SLDG-QC+BP)]{
\includegraphics[height=65mm]{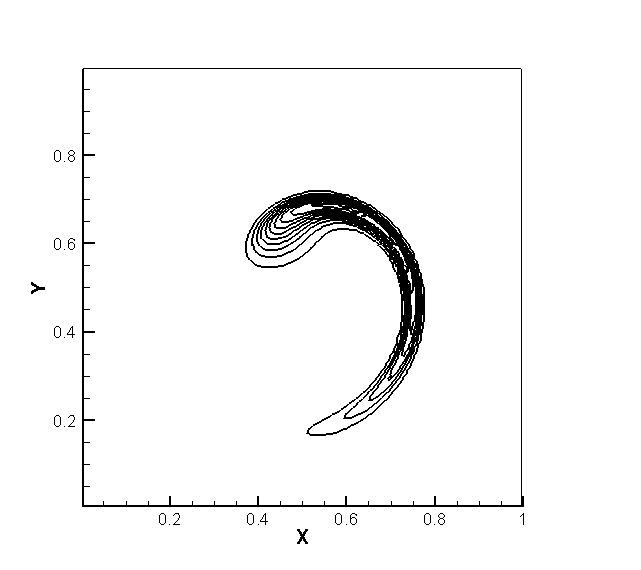} }
\subfigure[Field at T/2 ($P^2$ SLDG-QC+BP)]{
\includegraphics[height=65mm]{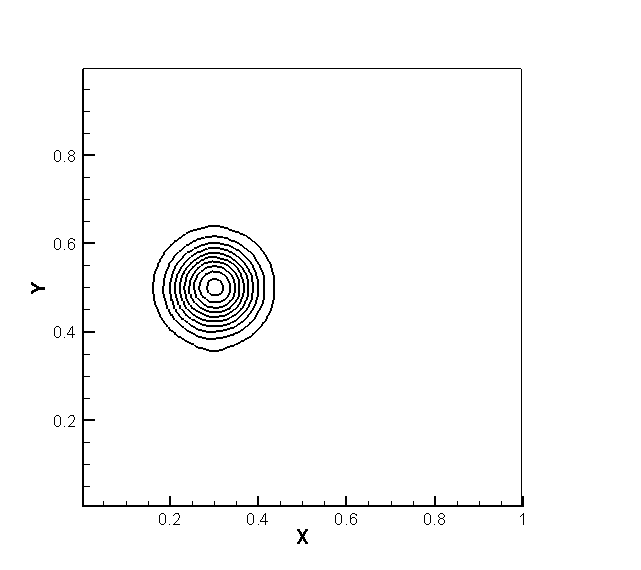} }
\subfigure[Field at 3T/4 ($P^2$ SLDG-QC+BP)]{
\includegraphics[height=65mm]{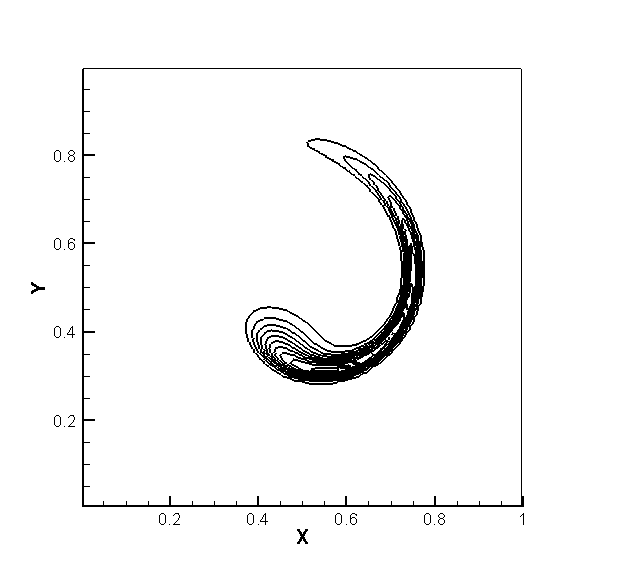} }
\caption{
2-D deformation flow. The $P^2$ SLDG scheme with quadratic-curved upstream cell with the BP filter is used.
The numerical mesh has a resolution of 60$\times$60. $\Delta t = 1.5\Delta x$.
Panels (b), (c), (d) and (a) shows the deformation of the initial distribution during the simulation at time $T/4$, $T/2$, $3T/4$ and $T$, respectively. The contours are plotted in the range from -0.05 to 0.95 with increment of 0.1, and an additional contour at 0 is added (see dashed line). The numerical solution is exactly positivity preserving. Thick contours are the highlighted exact (initial) solution for the contour values 0.05 and 0.75.
}
\label{df_p2QC_pp}
\end{figure}

\begin{table}[!ht]%\small
\caption{SLDG for 2D deformation flow. The
normalized $L^1$, $L^2$, $L^\infty$ errors and CPU time of SLDG schemes with or without the BP filter are reported at $T=2$ for comparison. $\Delta t=1.5\Delta x$. The numerical mesh has a resolution of 60$\times$60.}
\vspace{0.1in}
\centering
\begin{tabular}{l l ll c}

\hline
scheme & $L^1$ error & $L^2$ error & $L^\infty$ error  & CPU time (s)\\
\hline
 $P^1$ SLDG  & 4.03E-03 &      1.81E-02 &      2.23E-01  &  \textcolor[rgb]{1.00,0.00,0.00}{1.40} \\
 $P^1$ SLDG+BP  &  4.28E-03 &     2.05E-02 &    2.75E-01 & \textcolor[rgb]{1.00,0.00,0.00}{1.45} \\
 $P^2$ SLDG  &  7.41E-04 &      4.03E-03 &      1.43E-01 & \textcolor[rgb]{1.00,0.00,0.00}{3.15}      \\
 $P^2$ SLDG+BP  & 7.69E-04 &     4.21E-03 &     1.46E-01 & \textcolor[rgb]{1.00,0.00,0.00}{3.37}  \\
 $P^2$ SLDG-QC   & 4.82E-04 &  2.40E-03 &    3.24E-02  & \textcolor[rgb]{1.00,0.00,0.00}{3.82}      \\
 $P^2$ SLDG-QC+BP  & 5.24E-04 &    2.73E-03 &     4.55E-02 &   \textcolor[rgb]{1.00,0.00,0.00}{3.85}   \\
\hline
\end{tabular}
\label{table:rd_order}
\end{table}

\end{exa}

\section{Conclusion}
In this paper, we developed high order conservative semi-Lagrangian discontinuous Galerkin (SLDG) methods for 2D transport problems. The main ingredients of the proposed methods include the characteristic Galerkin weak formulation, the use of quadrilateral and the quadratic-curved quadrilateral approximations for the shape of upstream cells, and Green's theorem. The resulting schemes are up to third order accurate, locally conservative, unconditionally stable for linear constant coefficient problems, and are free of splitting errors in multi-dimensions. Accuracy and efficiency of the methods is showcased by classical tests for transport problems. Future work includes systematic theoretical investigation of the scheme for transport equations with variable coefficients, and development of the schemes for the Vlasov equation for plasma application and the multi-tracer transport model on the sphere.

\bigskip
\noindent\textbf{Acknowledgment}\\
We are grateful to Dr. Ram Nair from National Center for Atmospheric Research for helpful discussions.

%\appendix
%\input{append}

\bibliographystyle{abbrv}
\bibliography{references}

\end{document}